\documentclass[10pt,reqno,final]{amsart}
\usepackage{epsfig,amssymb,amsmath,version}
\usepackage{amssymb,version,graphicx,fancybox,mathrsfs}
\usepackage[notcite,notref]{showkeys}
\usepackage{url,hyperref}
\usepackage{subfigure}
\usepackage{color}
\usepackage{stmaryrd}
\usepackage{multirow}
\usepackage{booktabs,siunitx}
\usepackage{booktabs}
\usepackage{multicol}
\usepackage{booktabs} 
\usepackage{amsmath}
\usepackage{setspace}
\usepackage{array}
\usepackage{algorithm}
\usepackage{algorithmic}
\usepackage{enumerate}
\usepackage{enumitem}
\usepackage{tabu}                     
\usepackage{multirow}                 
\usepackage{multicol}                 
\usepackage{multirow}                
\usepackage{float}                    
\usepackage{makecell}                 
\usepackage{booktabs}                 
\usepackage{diagbox,bm}

\catcode`\@=11 \theoremstyle{plain}
\@addtoreset{equation}{section}   

\@addtoreset{figure}{section}
\renewcommand\thefigure{\thesection.\@arabic\c@figure}
\newtheorem{thm}{\bf Theorem}

\newtheorem{cor}{\bf Corollary}

\newtheorem{lmm}{\bf Lemma}

\theoremstyle{remark}
\newtheorem{rem}{\bf Remark}[section]

\definecolor{ligreen}{rgb}{0.0, 0.3, 0.0}

\definecolor{darkblue}{rgb}{0.0, 0.0, 0.55}

\definecolor{anti-flashwhite}{rgb}{0.55, 0.57, 0.68}

\newcommand{\bs}[1]{\boldsymbol{#1}}

\graphicspath{{./Figures/} } 

\begin{document}

\title[An Adaptive Orthogonal Basis Method for multiple solutions]{An Adaptive Orthogonal Basis Method for Computing Multiple Solutions of Differential Equations with polynomial nonlinearities}
\author[
	L. Li,\,    Y. Ye\,  $\&$\,  H. Li
	]{
		\;\; Lin Li${}^{\dag, *}$,   \;\;  Yangyi Ye${}^{\sharp}$ \;\; and \;\; Huiyuan Li${}^{\S}$
		}
	\thanks{${}^{\dag}$Corresponding author. School of Mathematics and Physics, University of South China, Hengyang, China. Email: exam-19861208@163.com (L. Li).\\
   \indent ${}^{*}$Key Laboratory of Computing and Stochastic Mathematics (Ministry of Education), School of Mathematics and Statistics, Hunan Normal University, Changsha, Hunan 410081, China.\\
   \indent ${}^{\sharp}$School of Mathematics and Physics, University of South China, Hengyang, China. Email: yeyangyi0911@163.com (Y. Ye)\\
		\indent ${}^{\S}$State Key Laboratory of Computer Science/Laboratory of Parallel Computing, Institute of Software, Chinese Academy of Sciences, Beijing, China. Email: huiyuan@iscas.ac.cn}



\begin{abstract} This paper presents an innovative approach, the Adaptive Orthogonal Basis Method, tailored for computing multiple solutions to differential equations characterized by polynomial nonlinearities. Departing from conventional practices of predefining candidate basis pools, our novel method adaptively computes bases, considering the equation's nature and structural characteristics of the solution. It further leverages companion matrix techniques to generate initial guesses for subsequent computations. Thus this approach not only yields numerous initial guesses for solving such equations but also adapts orthogonal basis functions to effectively address discretized nonlinear systems. Through a series of numerical experiments, this paper demonstrates the method's effectiveness and robustness. By reducing computational costs in various applications, this novel approach opens new avenues for uncovering multiple solutions to differential equations with polynomial nonlinearities.
\end{abstract}

\maketitle

\section{Introduction}\label{sect1}

In numerous mathematical models, particularly those involving nonlinear differential equations derived from real-world problems, the presence of nontrivial multiple solutions is a common occurrence. These multiple solutions frequently have direct relevance to practical applications \cite{2002A, davis1960introduction, Tadmor2012A}. However, it is widely acknowledged that providing explicit solutions in such cases is exceedingly challenging. As a result, researchers from around the globe often opt for the pursuit of numerical solutions. Therefore, the advancement and investigation of efficient numerical methods for the computation of multiple solutions take on paramount significance.

To the best of our knowledge, algorithms for computing multiple solutions can be broadly categorized based on the presence of a variational structure. Specifically, when dealing with differential equations possessing multiple solutions, the associated variational structure plays a pivotal role in shaping the algorithm's design for uncovering these multiple solutions. In 1993, Choi and McKenna introduced the mountain pass algorithm (MPA) for addressing multiple solutions in semilinear elliptic problems, drawing upon the mountain pass lemma in functional analysis \cite{choi1993mountain}. Subsequently, Xie et al. \cite{xie2005improved} highlighted the MPA's applicability in locating two solutions of mountain pass type, characterized by a Morse index of 1 or 0. In a different vein, Ding et al. underscored the MPA's limitation in computing sign-changing solutions and introduced a high linking algorithm (HLA) tailored to address such cases \cite{ding1999high}. Building on the foundational work of Choi, Ding, and others, in 2001, Zhou et al. proposed a local minimax method (LMM) inspired by the concepts presented in \cite{choi1993mountain, ding1999high}. The LMM characterizes a saddle point as a solution to a local minimax problem, offering another valuable approach to tackle multiple solution scenarios \cite{li2001minimax}. To be specific, Let $J(u)$ be a generic energy functional of differential equations with multiple solutions, and $J(u)$ is a $C^{1}$-functional on a Hilbert space $H$. Here it is worth pointing out that the solutions to differential equations with multiple solutions correspond to critical points of $J(u)$, and there exist saddle points belonging to critical points, where if $u^{*}$ is a saddle point of $J$, then we have
\begin{equation*}
  J(v) < J(u^{*}) < J(w), \quad \forall v, w \in N(u^{*}, \delta):=\{u|\; \|u - u^{*}\|\leq \delta \}\in H.
\end{equation*}
Based on the Morse index (MI) in the Morse theory, the LMM can obtain a saddle point with MI = $n$ ($n \in N^{+}$) by considering a two-level local minimax problem as follow:
\begin{equation}\label{eq1.3}
  \min_{v\in S_{H}}\max_{u\in[L, v]}J(u),
\end{equation}
where $S_{H}:=\{v\in H| \; \|v\| = 1 \}$ is the unit sphere. $L \subset H$ is a given ($n-1$) dimensional closed subspace, which can be constructed by some known critical points (or multiple solutions). $[L, v] := \{tv + w|\; t\geq 0,\; w\in L\}$ represents a closed half subspace. Some numerical algorithms can be conveniently implemented to solve \eqref{eq1.3}. In the LMM, the steepest descent direction is chosen as
the search direction in the local minimization process, and more recent advancements are presented in \cite{2005A, 2007numerical, 2008numerical, 2005Instability}. In cases where numerous differential equations exhibit multiple solutions without a discernible variational structure, the aforementioned methods become inapplicable. This circumstance gives rise to the second category of existing techniques for computing multiple solutions. The general procedure involves selecting certain numerical methods, such as the spectral method or finite difference method, to discretize the differential equations with multiple solutions. Subsequently, iterative methods are employed to locate the multiple solutions of the resulting nonlinear algebraic system (NLAS). In 2004, Xie at al. \cite{chen2004search} proposed the search-extension method (SEM), where they mainly considered the following nonlinear elliptic equations
\begin{equation}\label{1.2}
\begin{cases}
 -\Delta u + f(x, u) = 0,  \quad x\in \Omega,\\
  u = 0,  \quad\quad  \textrm{on}\; \partial\Omega,
  \end{cases}
\end{equation}
where $\Omega$ is a bounded domain in $R^{n}$ with a corresponding boundary $\partial \Omega$. In the search-extension method, the following eigenvalue problem is firstly solved
\begin{equation}\label{1.3}
\begin{cases}
 -\Delta \phi_{j} = \lambda_{j}\phi_{j},  \quad \textrm{in}\,\Omega,\\
  \phi_{j} = 0,  \quad\quad  \textrm{on}\; \partial\Omega,
  \end{cases}
\end{equation}
where $\{\lambda_{j}, \phi_{j}\} (j = 1, 2, \cdots)$ are its eigenpairs. Then the solution of \eqref{1.2} can be approximated by the following series:
\begin{equation}\label{1.4}
  u(x) = \sum_{j=1}^{N}a_{j}\phi_{j}.
\end{equation}
Substituting \eqref{1.4} into \eqref{1.2} yields the NLAS, and it is solved by the Newton method. Obviously, we can observe that $\{\phi_{j}\}_{j=1}^{\infty}$ from \eqref{1.3}-\eqref{1.4} are chosen to provide a good initial approximation of $u(x)$ in \eqref{1.2}. However, here it is worth pointing out that if the nonlinear term $f(x, u)$ in \eqref{1.2} plays a leading role, the choice of $\{\phi_{j}\}_{j=1}^{\infty}$ from \eqref{1.3}-\eqref{1.4} is not very suitable for a good initial approximation $u(x)$ in \eqref{1.2}.
As widely acknowledged, the Newton method suffers from a significant drawback, namely, its sensitivity to the initial guess and the conditioning of the Jacobian matrix. A recent study (\cite{2022Two}) introduced a promising alternative: the trust-region method, effectively replacing the Newton method for computing multiple solutions. This substitution not only significantly improved computational efficiency but also successfully addressed the aforementioned issues.

Moreover, the deflation technique has been utilized for the computation of multiple solutions (\cite{farrell2015deflation}). It is worth noting that in this context, the deflation procedure may encounter divergence problems even with a consistent initial guess. Additionally, an alternative approach involves using a randomly generated deviation from an already obtained solution as the initial guess to locate other solutions. However, it's crucial to highlight that these methods for initializing the guess may not be the most suitable choices, primarily because they do not adequately account for the underlying nonlinear characteristics inherent in the problem.

In addition, various other discretization approaches, such as finite difference methods, reduced basis methods, and finite element methods, have been coupled with homotopy continuation methods for computing multiple solutions \cite{allgower2006solution, allgower2009application, hao2014bootstrapping, 2018Two, zhang2013eigenfunction}. While it is true that the computational complexity escalates with mesh densification, ensuring the discovery of all solutions is of significant benefit. To address this challenge, a homotopy method with adaptive basis selection has been introduced \cite{hao2020homotopy}. Furthermore, in an effort to reduce computational costs in complex fields, companion matrix techniques have been leveraged for generating initial guesses \cite{hao2023companion}. Meanwhile, within the framework of a constructed dynamic system utilizing virtual time, the gentlest ascent dynamic (GAD) \cite{WE2011} and the constrained gentlest ascent dynamic (CGAD) \cite{WeiLiu2022} were proposed for computing multiple solutions. Subsequently, drawing inspiration from the shrinking dimer dynamic (SDD) \cite{2012JYZhang}, Zhang et al. introduced the (high-index) optimization-based shrinking dimer approach \cite{2016LeiZhang, 2019LeiZhang}. Additionally, the bifurcation method \cite{2007ZHYang, 2011ZXLi} was proposed for computing multiple solutions, drawing on principles from bifurcation theory.

In this paper, we mainly consider the following general differential equation with the nonlinearity of polynomial type, i.e.,
\begin{equation}\label{1.1}
\mathcal{L}{\bs{u}} + \mathcal{N}{\bs u} =
{\bs 0},  \quad\quad  {\bs x}\in \Omega,
\end{equation}
supplemented with some boundary conditions on $\partial \Omega$, where $\Omega \in \mathbb{R}^{d}$ is an open bounded domain, ${\bs u}:= (u_{1}({\bs x}), u_{2}({\bs x}), \cdots, u_{n}({\bs x}))^{T}$ is a vector function of ${\bs x}$ and $\mathcal{L}$ and $\mathcal{N}$ are linear and nonlinear operator, respectively. Here the linear operator $\mathcal{L}$ maybe $-\Delta$ or other operators. The nonlinear operator $\mathcal{N} := (\mathcal{N}_{1}({\bs u}), \mathcal{N}_{2}({\bs u}), \cdots, \mathcal{N}_{n}({\bs u}))^{T}$ is defined with nonlinearity of polynomial type. In the current work, we mainly focus on the case $n, d = 1$ or $2$. Traditionally, in spectral methods, the selection of basis functions may not always be well-suited for a given problem. Even though adaptive basis selection can help mitigate computational costs, it still entails choosing from a potentially extensive pool of candidate bases, as discussed in \cite{hao2020homotopy}. The core idea of this approach is to dynamically select the basis with the maximum residual based on the current solution using a greedy algorithm. Thus this technique effectively constructs a spectral approximation space tailored for nonlinear differential equations, and multiple solutions are then computed using the homotopy continuation method within this lower-dimensional approximation space.

In this paper, we expand upon the adaptive basis selection approach by introducing a novel adaptive basis method that deviates from the conventional practice of predefining a candidate basis pool. Instead, we dynamically compute the basis, taking into account both the nature of the equation and the structural characteristics of the solution. More specifically, we consider the polynomial nonlinearities within the differential equation and dynamically compute the basis functions to approximate multiple solutions. Once the basis is computed, we leverage companion matrix techniques, drawing inspiration from the work presented in \cite{hao2023companion}, to generate initial guesses for subsequent computation steps. This approach proves to be more efficient than the homotopy tracking method introduced in \cite{hao2020homotopy}.
Furthermore, we utilize the spectral trust-region method as a nonlinear solver for solving the NLAS \cite{2022Two}.

The structure of this paper is as follows: In Section \ref{sect2}, we provide an introduction to the fundamental concepts of the spectral collocation method and the trust-region method. Section \ref{sect3} introduces a new and innovative algorithm for computing multiple solutions. In Section \ref{sect4}, we present a comprehensive set of numerical experiments, discussing the efficiency and accuracy of our algorithm. Finally, we conclude the paper in Section \ref{sect5} with some closing remarks.

\section{Preliminary}\label{sect2}
To enhance the clarity and structure of our algorithm description in Section 3, it is essential to introduce some fundamental concepts. This section is divided into two parts: the first part explains the Spectral Chebyshev-Collocation method, and the second part outlines the trust-region method for iterating the resulting NLAS.

\subsection{The Spectral Chebyshev-Collocation method.}\label{sec2.1} The main idea of the Spectral Chebyshev-Collocation method is to use Lagrange interpolation based on Chebyshev points to approximate the solution of \eqref{1.1}. To account for boundary conditions, the Chebyshev-Gauss-Lobatto points $\{x_{j}\}_{j=0}^{N}$ are often chosen. These points $\{x_{j}\}_{j=0}^{N}$ are zeros of $(1-x^2)J^{1/2, 1/2}_{N-1}(x)$ where $J^{1/2, 1/2}_{N-1}(x)$ represents the Jacobi polynomial \cite{shen2011spectral}. This choice of Chebyshev points results in the following equation:
\begin{equation}
x_{j} = \cos\frac{j\pi}{N},  \quad  0\leq j \leq N.
\end{equation}
when we denote the values of the approximated function $u(x)$ at $\{x_{j}\}_{j=0}^{N}$ as
\begin{equation*}
  {\bs u} = (u(x_{0}), u(x_{1}), \cdots, u(x_{N}))^{T},
\end{equation*}
 the m-th derivative of $u(x)$ at $\{x_{j}\}_{j=0}^{N}$ (denoted by ${\bs u}^{(m)}$) can be expressed as
\begin{equation*}
  {\bs u}^{(m)} = D\cdot D\cdots D {\bs u} = D^{m}{\bs u},
\end{equation*}
Here, $D$ is known as the $(N+1)\times(N+1)$ Chebyshev spectral differentiation matrix, and its elements are as follows:
\begin{equation}\label{eq2.2}
\begin{split}
  ~ & (D)_{00} = \frac{2N^2+1}{6}, \quad\; (D)_{NN} = -\frac{2N^2+1}{6}, \\
   ~ & (D)_{jj} = \frac{-x_{j}}{2(1-x^2_{j})},  \quad  j = 1, \cdots, N-1,\\
   ~  & (D)_{ij} = \frac{c_{i}}{c_{j}}\frac{(-1)^{i+j}}{x_{i} - x_{j}},  \quad  i\neq j,  \quad  i, j = 1, \cdots, N-1,
\end{split}
\end{equation}
where
\begin{equation*}
c_{i} = \begin{cases}
2,   \quad\quad i=0\;\textrm{or}\; N,\\
1,  \quad\quad \textrm{otherwise}.
\end{cases}
\end{equation*}
For the case when $m = 2$, $D^{(2)} := D^{2} = D \cdot D$ can be derived, and its elements are as follows:
\begin{equation}\label{eq2.3}
(D^{(2)})_{ij} = \begin{cases}\vspace{0.1cm}
\frac{N^4 - 1}{15},  \quad \quad  & i = j = 0, \; i = j = N,\\\vspace{0.1cm}
-\frac{(N^2 - 1)(1 - x^2_{i}) + 3}{3(1-x^2_{i})^2},  \quad\quad & 1\leq i = j \leq N-1,\\\vspace{0.1cm}
\frac{2(-1)^{j}}{3c_j}\frac{(2N^2+1)(1-x_{j}) - 6}{(1-x_{j})^2},  \quad\quad   & i=0, \;  1\leq j \leq N,\\\vspace{0.1cm}
\frac{2(-1)^{j+N}}{3c_j}\frac{(2N^2+1)(1-x_{j})-6}{(1+x_{j})^2}, \quad\quad   & i = N, \; 0\leq j \leq N-1,\\
\frac{(-1)^{i+j}}{c_j}\frac{x^2_{i} + x_{i}x_{j} - 2}{(1-x^2_{i})(x_{i} - x_{j})^2}, \quad\quad  & 1 \leq i \leq N-1, \; 0\leq j \leq N, \; i\neq j.
\end{cases}
\end{equation}
Additionally, when differentiation matrices (e.g., \eqref{eq2.2} or \eqref{eq2.3}) are used for \eqref{1.1}, they must be combined with the corresponding boundary conditions (e.g., nonhomogeneous Dirichlet boundary conditions, Neumann boundary conditions). For simplicity, we have not provided these details here, and the reader is referred to \cite{Trefethenspectral} for more information.

To solve two-dimensional differential equations with multiple solutions, we extend the Chebyshev spectral differentiation matrices \eqref{eq2.2} - \eqref{eq2.3} to two-dimensional cases. As an illustrative example, we primarily focus on the Laplace operator \eqref{eq2.4}, namely,
\begin{equation}\label{eq2.4}
\Delta = \frac{\partial^2}{\partial x^2} + \frac{\partial^2}{\partial y^2}.
\end{equation}
Indeed, the Kronecker product, a technique from linear algebra, offers a convenient way to represent \eqref{eq2.4}. Specifically, for the second derivative with respect to $x$ in \eqref{eq2.4}, the corresponding Kronecker product is given by $D^{(2)}\otimes I$, where $I$ represents the identity matrix. This operation can be computed using the \textbf{kron}($D^{(2)}$, $I$) command in Matlab. Similarly, for the second derivative with respect to $y$ in \eqref{eq2.4}, the Kronecker product becomes $I \otimes D^{(2)}$. As a result, the discrete Laplace operator $L_{N}$ can be expressed as the Kronecker sum:
\begin{equation*}
  L_{N} = D^{(2)}\otimes I + I \otimes D^{(2)}.
\end{equation*}
This representation allows for the efficient computation of the Laplace operator in two-dimensional cases.

\subsection{The trust region method for iterating nonlinear system}\label{sec2.2}

Following the numerical discretization, such as the spectral Chebyshev-Collocation method discussed in Section \ref{sec2.1}, we obtain a nonlinear algebraic system represented  as:
\begin{equation}\label{eq2.2.4.1}
{\bs f} =\big({f}_{1}, \; {f}_{2}, \; \cdots, \; {f}_{n}\big)^{T}.
\end{equation}
Here, the unknown vector ${\bs a} = (a_{1}, a_{2}, \ldots, a_{n})^{T}$, with $\{f_i\}_{i=1}^{n}$ being functions of ${\bs a}$. Our objective is to find ${\bs a}$ that satisfies ${\bs f}({\bs a}) \equiv {\bs 0}$ using an iterative method. We can rewrite the system of nonlinear equations to a minimization problem below:
\begin{equation}\label{eq2.2.4.2}
\min_{\bs a \in {\mathbb R}^{n}}Q(\textbf{\textit{a}}),\quad Q(\bs a):= \frac{1}{2}\big\|\bs{f}(\bs a)\big\|^{2}_{2} = \frac{1}{2}\sum_{i = 1}^{n}{f}^2_{i}(\textbf{\textit{a}}).
\end{equation}
Consequently we have that for any ${\bs a}$, ${\bs f}({\bs a})\equiv 0$ iff $\min Q(\textbf{\textit{a}}) \equiv 0$. In other words, we can solve \eqref{eq2.2.4.2} to replace \eqref{eq2.2.4.1}.
First, we introduce the Jacobian matrix, the gradient, and the Hessian matrix, which are defined as follows:
\begin{equation*}\label{New4.4}
\bs J(\bs a) = {\bs f}'(\textbf{\textit{a}}) = (\nabla {f}_{1}(\textbf{\textit{a}}), \nabla {f}_{2}(\textbf{\textit{a}}), \cdots, \nabla f_{n}(\textbf{\textit{a}}))^{T},
\end{equation*}
\begin{equation*}
{\bs g}({\bs a}) = \nabla Q({\bs a}) = \bs{J}^{T}(\bs{a}) {\bs f(\bs {a})}   \quad \textrm{and} \quad   {\bs G}({\bs a}) = \nabla^2 Q({\bs a}) = {\bs J}^{T}(\bs{a}) \bs J({\bs a}) + \bs S({\bs a}),
\end{equation*}
where ${\bs S(\textbf{\textit{a}})} = \sum_{i = 1}^{n}f_{i}(\textbf{\textit{a}})\nabla^2 f_{i}(\textbf{\textit{a}})$.

The trust region method is based on the concept of defining a region around the current iteration where we trust the constructed quadratic model to be a suitable representation of the objective function $Q(\textbf{\textit{a}})$ in \eqref{eq2.2.4.2}. Then, we select an appropriate step to approximate the minimizer of the model within this region. In this method, both the direction and length of the step are considered simultaneously. Typically, the direction of the step changes whenever the size of the trust region is altered. If a step is deemed unacceptable, the size of the trust region is reduced, and a new minimizer is sought. The choice of the trust region size plays a pivotal role in the effectiveness of each step. If the trust region is too large, the minimizer of the constructed quadratic model may be far from the minimizer of the objective function $Q(\textbf{\textit{a}})$ in \eqref{eq2.2.4.2}. Conversely, if it's too small, the trust region method might miss an opportunity to take a substantial step that could bring it much closer to the minimizer of the objective function. Thus, it's essential to adjust the size of the region and retry. In practical computations, the size of the region is often adjusted based on the trust region method's performance in previous iterations. As explained in \cite{211Wright}, if the constructed quadratic model consistently provides reliable results, producing good steps and accurately predicting the behavior of the objective function along these steps, the size of the trust region can be increased to allow longer more ambitious steps. However, if the quadratic model is consistently unreliable, then the size of the trust region should be adjusted. After each such step, the size of the trust region should be reconsidered and potentially changed. The details of the trust region method for solving \eqref{eq2.2.4.2} are provided in the {\bf Appendix}.

\section{Adaptive Orthogonal Basis Method for computing multiple solutions of \eqref{1.1}}\label{sect3}

In this section we will focus on the main idea of the adaptive orthogonal basis method for computing multiple solutions of \eqref{1.1} and its general computational flow. Let's start with a set of orthogonal bases denoted as $\{\phi_i\}_{i=0}^{n-1}$, and the corresponding multiple solutions can be repressed as $u^{j}_{n-1} = \sum_{i=0}^{n-1}\alpha^{j}_i \phi_i$ ($j = 1, 2, \cdots$), where '$j$' represents the $j$-th solution in the current solution set and $\alpha^{j}_i$ represents the coefficient of the $j$-th solution. It is worth pointing out that these coefficients satisfy the discretized system $\bm f(\alpha_0^j,\cdots,\alpha_{n-1}^j) = \bm 0$ in \eqref{eq2.2.4.1}. To compute new solutions, a next basis function, denoted as $\phi_n$, is introduced to $\{\phi_i\}_{i=0}^{n-1}$, and we aim to ensure that new solutions can be expressed as $u_n(x) = \sum_{i=0}^n\alpha_i \phi_i$, with the condition that $\alpha_n\neq 0$, forcing the basis function $\phi_n$ to play a nontrivial role in computing new solutions. Note that the vector $\bm \alpha=(\alpha_0,\cdots,\alpha_{n})$ and the basis function $\phi_n$ are unknowns to be solved. To obtain them, an augmented system is constructed as follows:
\begin{equation}\label{AS}
\bm G(\bm \alpha,\phi_n)=
\begin{cases}
\frac{{\bm f}_n(\bm \alpha)}{\alpha_{n}} = {\bm 0}, \\
(\phi_{i},~ \phi_{n}) = 0, \quad i=0\cdots,n-1, \\
(\phi_{n},~ \phi_{n}) = 1.
\end{cases}
\end{equation}
In this context, our main aim is to a non-zero solution for
$\alpha_{n}$, ensuring that the newly introduced basis function $\phi_n$ significantly contributes to new solutions.

Next, we will compute multiple solutions of \eqref{AS}. due to the high nonlinearity and the multiplicity of the solution, solving multiple solutions of \eqref{AS} can be extremely challenging. In such cases, our method is to leverage existing solution information (denoted as $\alpha_0^j,\cdots,\alpha_{n-1}^j$), i.e. they are used as the initial guess to solve the unknown vector $\bm a$. Meanwhile,
the variable $\alpha_n$ can be obtained by solving a
polynomial algebraic equation, i.e. the roots of the polynomial algebraic equation are the values of the variable $\alpha_n$. This simplification has many advantages. On the one hand, it reduces the computational complexity. On the other hand, it is especially useful when dealing with polynomials whose roots can be determined using the eigenvalue method via a companion matrix. To be specific, we consider a
general monic polynomial of degree $n$ given by:
\begin{equation}
p(t) = a_0 + a_1t + \ldots + a_{n-1}t^{n-1} + t^n,
\end{equation}
where $\{a_{i}\}_{i=0}^{n-1}$ are real numbers. Consequently, the corresponding companion matrix, denoted as $C_{p}$, takes the form:
\begin{equation}
C_p =
\begin{bmatrix}
0 & 1 & 0 & \cdots & 0 \\
0 & 0 & 1 & \cdots & 0 \\
\vdots & \vdots & \vdots & \vdots & \vdots \\
0 & 0 & 0 & \cdots & 1 \\
-a_0 & -a_1 & -a_2 & \cdots & -a_{n-1}
\end{bmatrix}.
\end{equation}
Then, the eigenvalues of the companion matrix $C_{p}$ are the values of the variable $\alpha_n$, and this can be efficiently computed by \textbf{eig}($C_{p}$) in MATLAB.\vspace{0.2cm}

Now we state the main idea of the adaptive orthogonal basis method as follows:

\begin{enumerate}
\item \textbf{Initialization}: Generate Chebyshev-Gauss-Lobatto points on the domain $\Omega$. Next, we solve Eq. \eqref{eq2.2.4.2} by using the trust region method outlined in Section \ref{sec2.2}, and the resulting numerical solution is denoted as ${\bm u}_{0}$. One can then set $\alpha_{0} = \|{\bm u}_{0}\|_{2}$ and $\phi_{0}={\bm u}_{0}/\alpha_0$, which forms a unit basis function.

\item \textbf{Orthogonal Basis expansion iteration}: Consider a solution set denoted as $\mathcal{S} = {({\bm \alpha}^1, {\bm \alpha}^2, \cdots)}$. Each element in $\mathcal{S}$ corresponds to
    a solution of \eqref{1.1} with a set of basis functions $(\phi_0, \ldots, \phi_{n-1})$, i.e. $u^{j} = \sum_{i=1}^{n-1}\alpha^j_i\phi_i$. To obtain new solutions, the orthogonal basis function $\phi_n$ should be introduced into $\{\phi_i\}_{i=1}^{n-1}$. As a result, these new solutions can be expressed as $\hat{u}^{j} = \sum_{i=0}^{n}\alpha^j_i\phi_i$. For each $j$-th solution, we compute $\alpha_n$ and $\phi_n$ by solving the augmented system described in Equation \eqref{AS} using the trust region method. Subsequently, for each $\phi_n$ obtained, we calculate multiple solutions for $\alpha_n$ by computing the eigenvalues of the companion matrix.

\item {\bf Local orthogonalization:} When we compute $\phi_n$, more than one basis function may be obtained due to multiple solutions, where these basis functions $(i.e. \phi_n, \ldots, \phi_{m})$ are called as adaptive basis functions, and are not orthogonal. The Gram-Schmidt orthogonalization method is used to generate a set of adaptive orthogonal basis functions. Afterward, we can expand a set of adaptive orthogonal basis functions $\{\phi_0, \ldots, \phi_{n-m}\}$.

\item {\bf Filtering conditions:} In the computational process, some spurious solutions may be obtained. Based on the collocation method, we use the residual to eliminate these spurious solutions, thereby reducing the computational cost.

\item {\bf Stopping criteria:} Our algorithm will terminate when no more basis can be computed in the augmented system \eqref{AS}.
\end{enumerate}

Next, we are ready to state a general computational flow for computing multiple solutions of \eqref{1.1} as follows:
\begin{itemize}
  \item[\textbf{Step 1}.] Based on the spectral-collocation method to \eqref{1.1}, a nonlinear algebraic system is obtained and shown in \eqref{eq2.2.4.1}. With the trust region method presented in section \ref{sect2}, the resulted numerical solution $u^{(0)} := \alpha_0 \phi_0$ can be obtained, where $\alpha_0 = \|u^{(0)}\|_2$, and $\phi_0$ is a unit basis function.
  \item[\textbf{Step 2}.] Let $\hat{u}^{(0)} := \alpha \phi_0$, and substitute it into \eqref{1.1}. a polynomial algebraic equation on $\alpha$ is obtained by using the spectral-Galerkin method. Based on the eigenvalues of the corresponding companion matrix $C_p$, multiple solutions on $\alpha$ can be obtained, and denoted by $\{\alpha^{(0)}_0, \alpha^{(1)}_0, \cdots\}$, which leads to $\hat{u}^{(0)}_i = \alpha^{(i)}_0 \phi_0\;(i = 0, 1, \cdots)$.
  \item[\textbf{Step 3}.] We set $\hat{u}^{(1)} = \alpha_{1, 0}\phi_0 + \alpha_{1, 1}\phi_1$ to obtain multiple solutions with high precision, where $\alpha_{1, 0}, \alpha_{1, 1}$ and $\phi_1$ are unknown variables to be solved. The augmented system \eqref{AS} is used to determine them, where the initial guess on $\alpha_{1, 0}$ is taken from $\{\alpha^{(0)}_0, \alpha^{(1)}_0, \cdots\}$. As a result, we can obtain multiple solutions, and denote them by $\{\alpha_{1, 0}^{(j)}, \alpha_{1, 1}^{(j)}, \phi^{(j)}_1\}\;(j = 1, 2, \cdots, k)$.
  \item[\textbf{Step 4}.] The Gram-Schmidt orthogonalization method is used to $\{\phi^{(j)}_1\}_{j = 1}^{k}$, which leads to a set of adaptive orthogonal basis functions, denoted by $\{\phi_1, \phi_2, \cdots, \phi_k \}$. Subsequently, we set the numerical solution $\tilde{u}^{(1)} = \alpha_{1, 0}\phi_0 + \tilde{\alpha}_{1, 1}\phi_1$, where $\tilde{\alpha}_{1, 1}$ is a single unknown to be solved. Similar to $\alpha$ in \textbf{Step 2}, $\tilde{\alpha}_{1, 1}$ is obtained.
  \item[\textbf{Step 5}.] To obtain multiple solutions with high precision, we will introduce $\{\phi_i \}_{i = 2}^{k}$ into $\tilde{u}^{(1)}$, i.e. $\tilde{u}^{(1)} = \alpha_{1, 0}\phi_0 + \tilde{\alpha}_{1, 1}\phi_1 + \sum_{j=2}^{l}\alpha_{1, j}\phi_{j} \; (2\leq l\leq k)$, and repeat from \textbf{Step 2} to \textbf{Step 4}. Multiple solutions to \eqref{1.1} can be obtained.
\end{itemize}

Compared with the existing methods, the main differences and advantages of our algorithm reside in the following aspects:

\begin{itemize}
    \item {\bf Enhanced Initial Guesses:} In general, in the nonlinear iteration, a good initial guess often plays a vital role in computing solutions for differential equations, particularly when dealing with multiple solutions. Our algorithm introduces a significant improvement by generating multiple appropriate initial guesses based on the eigenvalues of the companion matrix. To be specific, In \textbf{Step 2}, $\alpha$ can be obtained from the companion matrix $C_p$. Similar situation is also seen for $\tilde{\alpha}_{1, 1}$ in \textbf{Step 4}.
    This approach deviates from existing methods like the search extension method \cite{chen2004search} and enhances the reliability and diversity of a good initial guess.
\item {\bf Robustness with the Trust Region Method:}  Our approach is based on employing the trust region method when iterating through nonlinear algebraic systems. The trust region method allows for more relaxed and flexible choices of initial inputs. This feature proves invaluable as the complexity of the nonlinear algebraic system increases, ensuring robust performance even as more solutions are sought.

\item {\bf Reduced the Computational Complexity:} In our algorithm, we dynamically compute the orthogonal basis by solving the augmented system instead of relying on pre-defined basis functions (see \textbf{Step 4}). With a good initial guess, this adaptive approach reduces the basis function set, and doesn't increase the number of unknown variables, indicating that the computational complexity will be reduced.
\end{itemize}

In summary, our adaptive orthogonal basis algorithm not only provides improved initial guesses and computational efficiency but also leverages the trust region method for robust convergence. Additionally,
based on adaptive orthogonal basis functions, our algorithm can reduce the computational complexity. These advantages set our algorithm apart from existing algorithms and establish it as a valuable method for solving equations with multiple solutions.

\section{Numerical experiments}\label{sect4}
In this section, our primary objective is to assess and demonstrate the effectiveness of our algorithm. We present one-dimensional and two-dimensional examples in Sections \ref{sect4.1} and \ref{sect4.2}, respectively. All code execution is performed on a server equipped with an Intel(R) Core(TM) i7-11700F processor running at 2.50 GHz and 32GB of RAM, utilizing MATLAB (version R2020b).

\subsection{1D Examples.\vspace{0.25cm}\\}\label{sect4.1}

\emph{Example 1.} We consider the following nonlinear boundary value problem:
\begin{equation}\label{eq4.1}
\begin{cases}\vspace{1.5mm}
u^2_{xx} - (1+e^{x})u_{xx} + e^{x} = 0,  \quad\quad  x\in (0, 1),\\
u(0) = u(1) = 0.
\end{cases}
\end{equation}
Upon solving the algebraic differential equation, we obtain two linear differential equations: $u_{xx}=1$ and $u_{xx}=e^x$. Consequently, we have two corresponding solutions:
\begin{equation}\label{eq4.2}
u_{1} = \frac{x^2}{2}-\frac{x}{2} \hbox{~and~}  u_{2} = e^x-(e-1)x-1.
\end{equation}
Based on the adaptive orthogonal basis method presented in section \ref{sect3}, the corresponding flowchart is shown in Fig.\ref{Figure4.1}, and some further explanations are presented as follows:\vspace{0.2cm}
\begin{figure}[!ht]\label{Figure4.1}
	\centering
	\includegraphics[scale=.4]{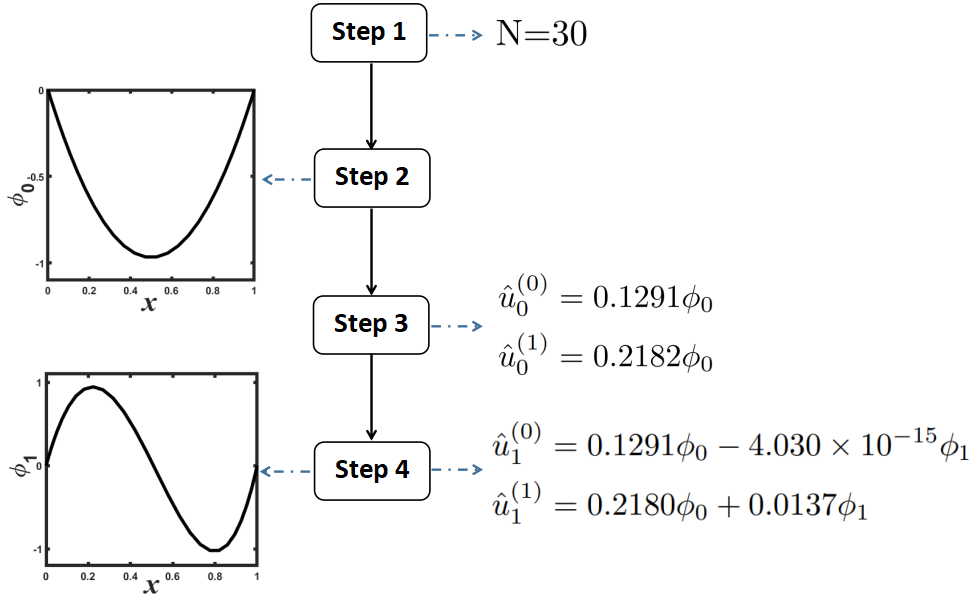}
	\caption{The flow chart of the adaptive orthogonal basis method for computing multiple solutions of  \eqref{eq4.1}. Two basis functions, $\phi_{0}$ and $\phi_{1}$, form two solutions, $\hat{u}^{(0)}_{1}$ and $\hat{u}^{(1)}_{1}$, respectively, in Step 4. A notable difference between $\hat{u}^{(0)}_{1}$ and $\hat{u}^{(1)}_{1}$ lies in the coefficients of the basis function $\phi_{1}$. Specifically, the coefficient of $\phi_{1}$ for $\hat{u}^{(0)}_{1}$ is considerably smaller than that for $\hat{u}^{(1)}_{1}$. This suggests that the second solution, $\hat{u}^{(1)}_{1}$, may be formed by extending the first solution space ${\phi_{0}}$ through the inclusion of an adaptive basis function, $\phi_{1}$.  }
\end{figure}

%
%
%
%
%

\begin{itemize}
  \item In Step 3, when applying the spectral Legendre-Galerkin method, it is essential to transform the function values initially computed at Chebyshev-Gauss-Lobatto points (as described in Step 2) into values at Legendre-Gauss-Lobatto points. The transformation process is depicted in Figure \ref{Fg3.1}. It is noteworthy that various fast, straightforward, and numerically stable algorithms for Chebyshev-Legendre transforms can be found in the literature, including references such as \cite{1974Ahmed, 1991Alpert, 1998Potts, 2014HaleandTownsend, shen2011spectral}.

\begin{figure}[!ht]
\begin{centering}
\includegraphics[width=14cm,height=1.8cm]{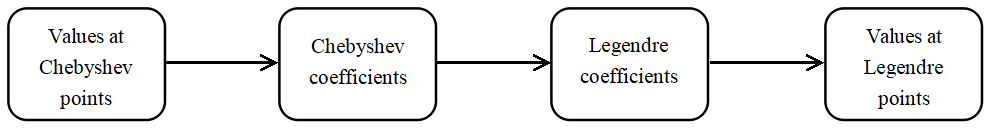}
\caption{\small Illustration of Chebyshev-Legendre transforms.}\label{Fg3.1}
\end{centering}
\end{figure}

  \item In Step 4, our primary objective is to enhance the accuracy of multiple solutions by increasing the number of standard orthogonal basis functions, such as $\phi_1$. In other words, a sequence of standard orthogonal basis functions ${\phi_i}$ with $i=1, 2, \cdots$ can be systematically increased until the predefined stopping criteria are satisfied. Additionally, the values $\alpha^{(i)}_0$ obtained in Step 3 can serve as initial guesses for $\alpha_{1, 0}$ in the process of solving \eqref{AS}.

\end{itemize}

\begin{figure}[ht]\label{final}
  \centering
  \subfigure[Numerical errors between exact and numerical solutions with N=30]{\includegraphics[width=5.8cm]{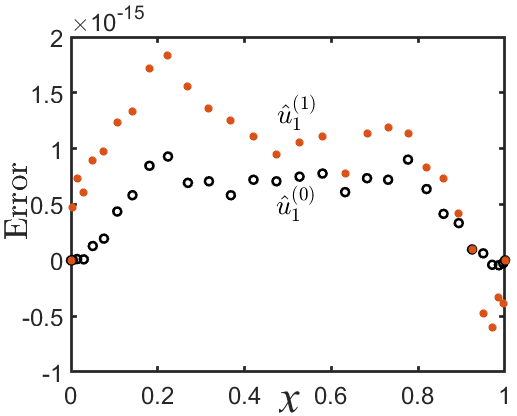}}\;\;\quad\quad\quad
  \subfigure[Convergence of two numerical solutions]{\includegraphics[width=7cm,height=4.5cm]{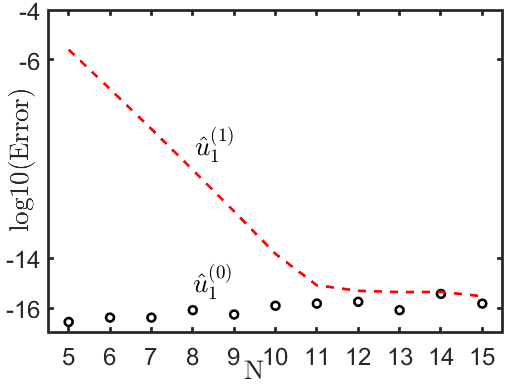}}\;\;
  \caption{Two solutions and its convergence for \eqref{eq4.1}.}\vspace{-0.8cm}
\end{figure}

In Fig.\ref{final}a, we plot the numerical errors for these two solutions.
In Fig.\ref{final}b, we present the convergence test for both solutions. The first solution, $\hat{u}^{(0)}_{1}$, demonstrates machine accuracy when $N$ is very small, while the second solution, $\hat{u}^{(1)}_{1}$, exhibits spectral accuracy. This distinction arises due to the notably low algebraic degree of the first solution, as evident in $u_{1}$ from \eqref{eq4.2}. Furthermore, when varying $N$ across values of 5, 10, 15, and 20, our algorithm for computing the solutions showcases remarkable computational efficiency, as demonstrated in Table \ref{Table1}.

\begin{table}[H]
	\centering
	\caption{Computational times for Example 1 with different $N$.}\vspace{0.4cm}
\begin{tabular}{|ccccc|}
  \hline
  $N$ & 5 & 10 & 15 & 20 \\\hline
  Time(s) &0.48 &0.52 & 0.53 &0.60 \\
  \hline
\end{tabular}\label{Table1}\vspace{-0.1cm}
\end{table}

%

\begin{figure}[ht]
  \centering
  \subfigure[The basis function $\phi_0$]{\includegraphics[width=5cm,height=4.1cm]{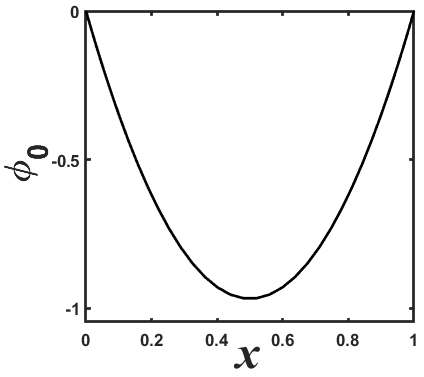}}\;\;\quad\quad\quad
  \subfigure[Numerical errors between exact and numerical  solutions with N=30]{\includegraphics[width=5.5cm]{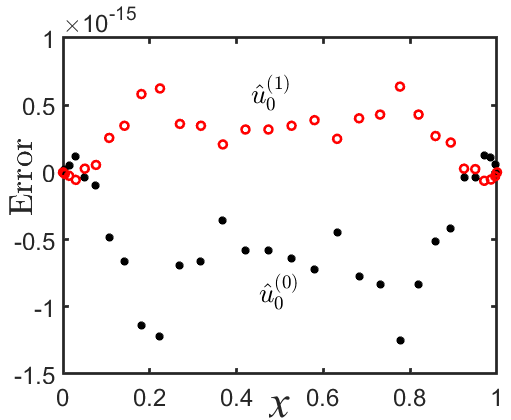}}\;\;
  \caption{Multiple solutions of  \eqref{eq4.3} by using our method.}\label{eq4.2figure2}
\end{figure}

\emph{Example 2.} The following boundary value problem is considered:
\begin{equation}\label{eq4.3}
\begin{cases}\vspace{0.15cm}
u^2_{xx}+u_{xx}-2=0 ,\;\;\;\;\quad  x\in (0,1),\\
u(0)=0,\;\;u(1)=0,
\end{cases}
\end{equation}
which has two solutions, i.e.,
\begin{equation}\label{eq4.4}
u^+=\frac{x^2}{2}-\frac{x}{2},\;\;\;\;u^-=-x^2+x.
\end{equation}
Obviously, we can conclude $\frac{u^{+}}{u^{-}}= -\frac{1}{2}$.
Employing our algorithm, we identify an adaptive basis function $\phi_{0}$, as shown in Fig. \ref{eq4.2figure2}a. Additionally, we determine two coefficients: $\alpha^{(0)}_0 = -0.2582$ and $\alpha^{(1)}_0 = 0.1291$, from which we form two numerical solutions $\hat{u}^{(i)}_{0} = \alpha^{(i)}_0\phi_{0}$ for $i = 0, 1$. A comparison between the exact and numerical solutions is presented in Fig. \ref{eq4.2figure2}b. In Table \ref{eq4.2figure3}, we provide details on numerical errors, residuals, and computational times for the two numerical solutions (labeled as I and II). These results affirm the feasibility and effectiveness of our algorithm.

\begin{table}[h]
	\centering
	\caption{Performance of our algorithm for solving Eq. \eqref{eq4.3}.}\label{eq4.2figure3}
	\begin{tabular}{|c|cc|cc|c|}
		\hline
 \multirow{2}{*}{N}&\multicolumn{2}{c}{Solution I }  &\multicolumn{2}{c|}{Solution II }  &\multirow{2}{*}{Time(s)}\\
 \cline{2-5}
		&Error &Residual &Error &Residual  &  \\\hline
5   & 5.5511e-17 & 2.6645e-15 &    2.7756e-17   &    1.3323e-15       &   0.3187 \\
    	10  & 3.3307e-16 & 1.7319e-14  & 1.6653e-16& 6.6613e-15   &   0.3310 \\
		15	& 2.4980e-16 & 6.5281e-14 &  1.5265e-16& 2.1316e-14   &   0.3511 \\\hline
	\end{tabular}
\end{table}

\begin{figure}[!h]
	\begin{center}
		\subfigure[$\lambda=1,\, p=18$ ]{ \includegraphics[scale=.26]{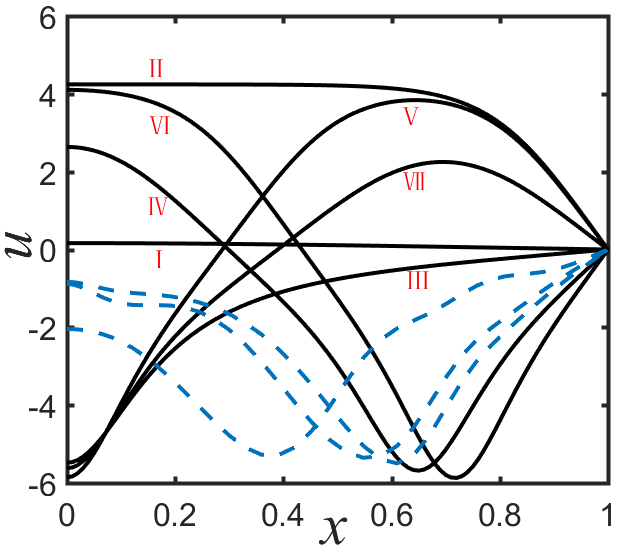}}\quad\quad\;\;
		\subfigure[$\lambda=-\pi^2/4,\, p=10$ ]{\includegraphics[scale=.26]{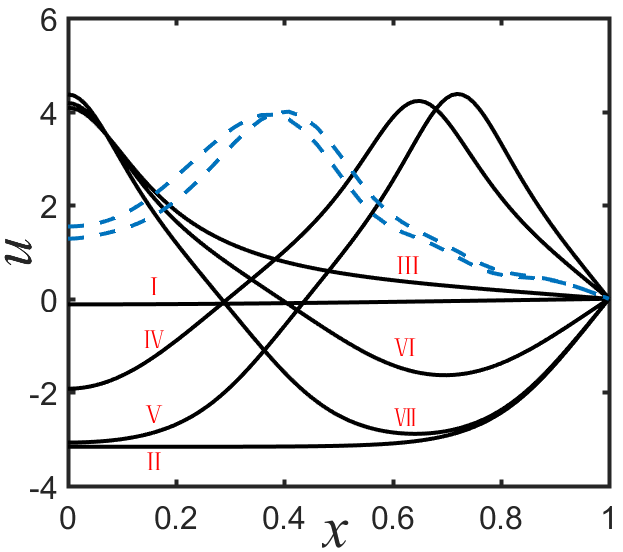}}
		\caption{\small Multiple solutions of \eqref{eq4.5} with different $(\lambda,\,p)$, and some spurious solutions (blue line) before applying filtering conditions.}\label{Newfigure1}
	\end{center}
\end{figure}

\begin{figure}[!h]
	\begin{center}
		\subfigure[$\phi_0$ ]{ \includegraphics[scale=.185]{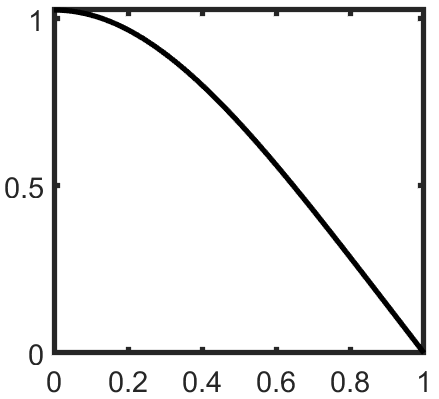}}\;\;
		\subfigure[$\phi_1$ ]{\includegraphics[scale=.185]{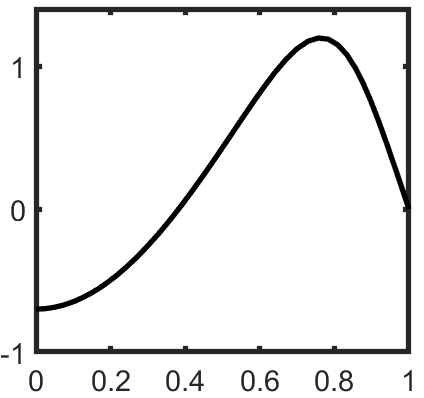}}\;\;
           \subfigure[$\phi_2$ ]{\includegraphics[scale=.185]{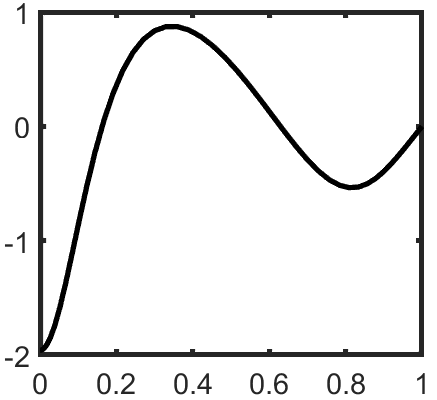}}\;\;
           \subfigure[$\phi_3$ ]{\includegraphics[scale=.185]{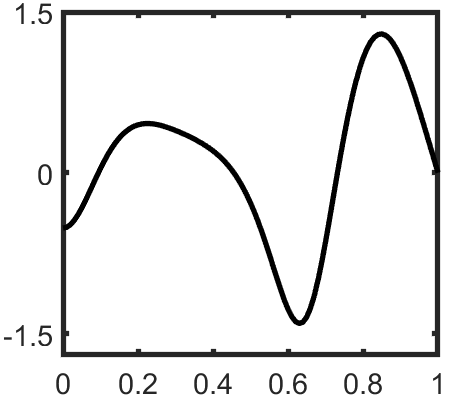}}\\
           \subfigure[$\phi_4$ ]{\includegraphics[scale=.185]{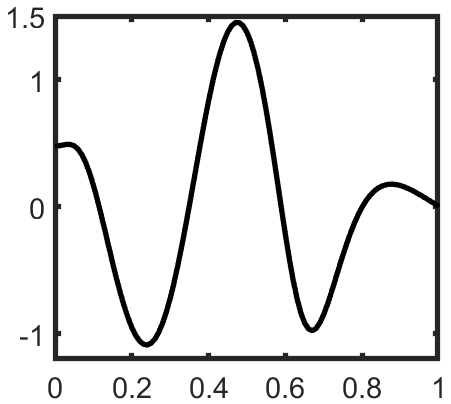}}\quad\quad
          \subfigure[$\phi_5$ ]{\includegraphics[scale=.185]{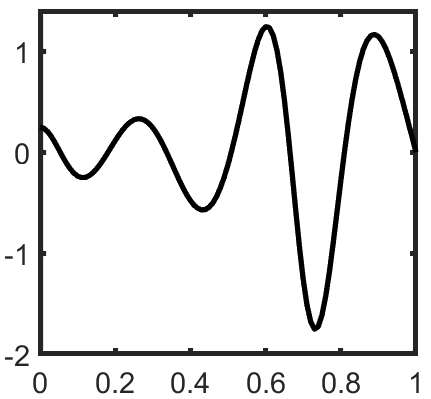}}\quad\quad
          \subfigure[$\phi_6$ ]{\includegraphics[scale=.185]{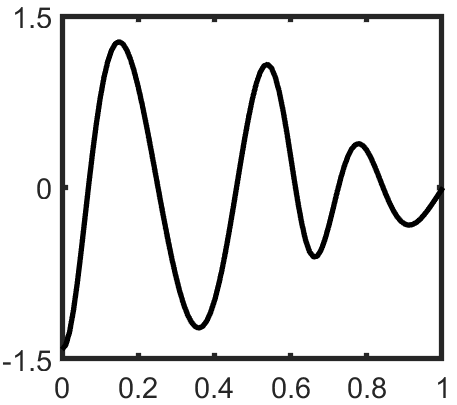}}
		\caption{\small Basis functions computed by our algorithm for solving \eqref{eq4.5} with $\lambda=1$ and $p=18$.}\label{Newfigure2}
	\end{center}
\end{figure}

\begin{figure}[!h]
	\begin{center}
		\subfigure[$\phi_0$ ]{ \includegraphics[scale=.185]{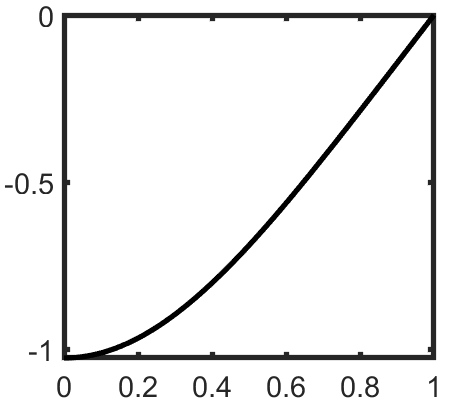}}\;\;
		\subfigure[$\phi_1$ ]{\includegraphics[scale=.185]{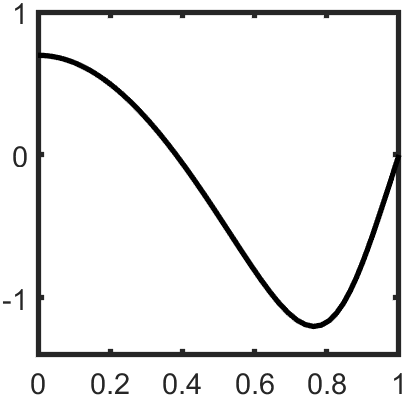}}\;\;
           \subfigure[$\phi_2$ ]{\includegraphics[scale=.185]{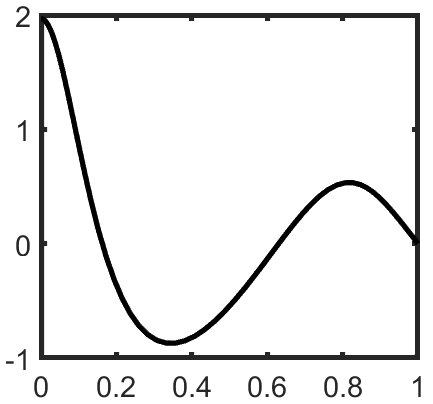}}\;\;
           \subfigure[$\phi_3$ ]{\includegraphics[scale=.185]{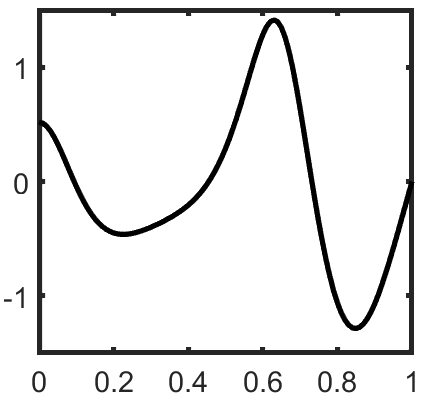}}\\
           \subfigure[$\phi_4$ ]{\includegraphics[scale=.185]{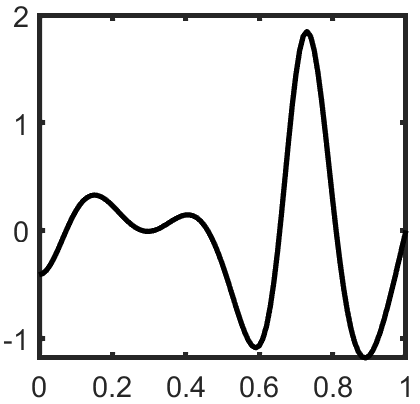}}\quad\quad
          \subfigure[$\phi_5$ ]{\includegraphics[scale=.185]{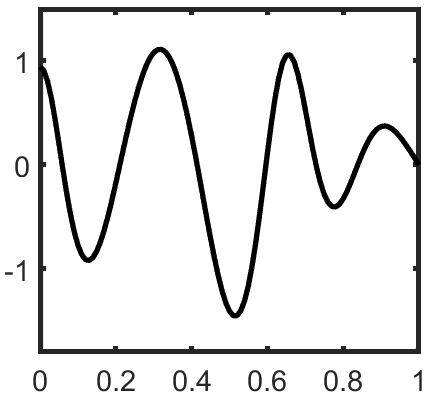}}\quad\quad
          \subfigure[$\phi_6$ ]{\includegraphics[scale=.185]{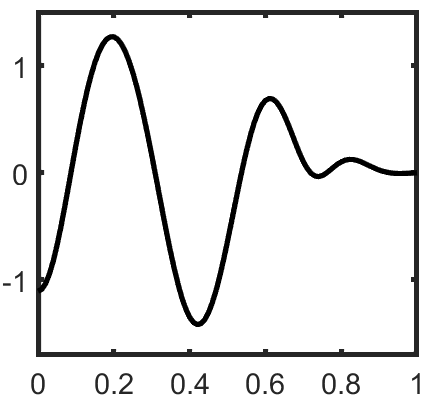}}
		\caption{\small Basis functions computed by our algorithm for solving \eqref{eq4.5} with  $\lambda=-\pi^2/4$ and $ p=10$.}\label{Newfigure3}
	\end{center}
\end{figure}

\emph{Example 3.} Next we consider the following boundary value problem with multiple solutions\cite{hao2023companion}:
\begin{equation}\label{eq4.5}
\begin{cases}\vspace{0.15cm}
u_{xx}=\lambda u^2(u^2-p) ,\;\;\;\;\quad  x\in (0,1),\\
u'(0) = 0, \quad\; u(1) = 0.
\end{cases}
\end{equation}
Here, the parameters $\lambda$ and $p$ are introduced. Employing our algorithm with $\lambda=1$ and $p=18$, we obtain multiple solutions, as depicted in Fig. \ref{Newfigure1}a. These numerical solutions $u_{N}(x)$ can be expressed using the following adaptive basis functions:
\begin{equation}\label{eq4.6}
u_{N}(x) = \sum_{i=0}^{6}\tilde{\alpha}_{6,i} \phi_{i}.
\end{equation}
The adaptive basis functions and coefficients ${\tilde{\alpha}_{6, i}}$ as defined in \eqref{eq4.6} are illustrated in Fig. \ref{Newfigure2} and tabulated in Table \ref{Newtable1}, where the coefficients above the diagonal are very small, indicating that the difference between these solutions can be reflected in the adaptive basis functions. In addition, our algorithm effectively excludes some spurious solutions by applying the filtering conditions, leading to excellent agreement with the solutions previously reported in \cite{hao2014bootstrapping}.
A similar observation can be made when setting $\lambda=-\pi^2/4$ and $p=10$, as shown in Fig. \ref{Newfigure3} and Table \ref{Newtable2}.

\begin{table}[!h]
	\renewcommand{\arraystretch}{1.2}  	\setlength\tabcolsep{1.3pt}
	\centering
	\caption{Numerical results of \eqref{eq4.5} with $\lambda=1,\;p=18$. }\label{Newtable1}
	\begin{tabular}{|c|ccccccc|c|c|}
		\hline
		\diagbox{Sol.ind}{Coefs} & $\tilde{\alpha}_{6,0}$ & $\tilde{\alpha}_{6,1}$ & $\tilde{\alpha}_{6,2}$ & $\tilde{\alpha}_{6,3}$ & $\tilde{\alpha}_{6,4}$ & $\tilde{\alpha}_{6,5}$ & $\tilde{\alpha}_{6,6}$ &Residual & Time(s) \\ \hline
		
		I & 0.16 & {\color{red}3.45e-15} & {\color{red}3.36e-15} & {\color{red}-4.05e-16} & {\color{red}3.87e-16} & {\color{red}-5.55e-17} & {\color{red}1.88e-16} & 2.14e-12  &\multirow{7}{*}{10.83}\\
		II &  5.17 & 1.51 &  {\color{red}1.07e-13} & {\color{red}1.83e-15} & {\color{red}-1.40e-14} & {\color{red}-9.70e-15} & {\color{red}3.07e-14} & 4.15e-11 &~\\
		III &  -2.57 & 0.93 &  1.12 & {\color{red}1.18e-12} &  {\color{red}-1.82e-12} & {\color{red}-4.73e-13} &  {\color{red}4.85e-12} & 1.32e-11 &~\\
		IV & -0.90  & -3.76 & -0.78 & 1.19 & {\color{red}-9.92e-06}  & {\color{red}-2.44e-06}  & {\color{red}-2.00e-05} & 4.10e-12 &~\\
		V & 0.08 & 3.85 &  1.83 &  -0.35 & 0.34 &  {\color{red}-4.70e-11} & {\color{red}-6.54e-11} & 2.13e-11 &~\\
		VI &  1.08 & -4.57 &  0.01 & 0.70 & 0.15 & 0.43 & {\color{red}-5.07e-06} & 1.55e-11  &~ \\
		VII & -1.28 & 2.72 & 1.26  & -0.24 & 0.04 & -0.01 & 0.06 & 8.34e-11 &~ \\
		\hline
	\end{tabular}
\end{table}

\begin{table}[!h]
	\renewcommand{\arraystretch}{1.2}
	\setlength\tabcolsep{1.3pt}
	\centering
	\caption{Numerical results of \eqref{eq4.5} with $\lambda=-\pi^2/4,\; p=10$. }\label{Newtable2}
	\begin{tabular}{|c|ccccccc|c|c|}
		\hline
		\diagbox{Sol.ind}{Coefs} & $\tilde{\alpha}_{6,0}$ & $\tilde{\alpha}_{6,1}$ & $\tilde{\alpha}_{6,2}$ & $\tilde{\alpha}_{6,3}$ & $\tilde{\alpha}_{6,4}$ & $\tilde{\alpha}_{6,5}$ & $\tilde{\alpha}_{6,6}$ &Residual & Time(s) \\ \hline
		I & 0.12 & {\color{red}8.43e-15} & {\color{red}7.24e-15} & {\color{red}-1.19e-15} & {\color{red}3.28e-16} & {\color{red}4.89e-16} & {\color{red}1.03e-15} & 8.59e-13 &\multirow{7}{*}{17.36}\\
		II &  3.86 & 1.14 &  {\color{red}5.78e-14} & {\color{red}-5.47e-15} & {\color{red}-4.23e-15} & {\color{red}7.49e-15} & {\color{red}-6.80e-15} &4.44e-11  &~\\
		III &  -1.90 & 0.69 &  0.84 & {\color{red}-3.07e-14} &  {\color{red}-2.79e-15} & {\color{red}-8.00e-14} &   {\color{red}-1.00e-13} & 1.16e-11 &~ \\
		IV & -0.69  & -2.77  & -0.59 & 0.90  & {\color{red}-1.44e-06}  & {\color{red}-1.04e-06} & {\color{red}-6.70e-07} & 2.89e-11 &~ \\
		V & 0.87 & -3.40 &  0.04 &  0.50 & 0.36 &  {\color{red}2.31e-07} & {\color{red}-3.50e-07} & 2.60e-11  &~\\
		VI &  -0.97 & 2.00 &  0.94 &- 0.18& 2.36e-04& 0.05 & {\color{red}1.21e-13}  & 7.67e-11   &~\\
		VII &  0.11 & 2.88 &  1.39 &- 0.27 & 0.08 & 0.14 & 0.21 & 3.04e-11 &~\\
		\hline
	\end{tabular}
\end{table}

\subsection{2D Examples.\vspace{0.25cm}\\}\label{sect4.2}

\emph{Example 4.} We first consider the following 2D example:
\begin{equation}\label{eq4.7}
\begin{cases}\vspace{0.15cm}
\Delta{u}+u^2=800\sin(\pi x)\sin(\pi y) \;\;\; &{\rm in} \;\;\;\Omega = (0,1)\times(0,1),\\
u=0  &{\rm on}\;\;\;\partial\Omega.
\end{cases}
\end{equation}
By employing our algorithm, we have successfully computed 10 solutions, as depicted in Fig. \ref{New2figure2}. This results are in line with the findings of Breuer et al. in their work \cite{2003Breuer}, where it was established that Eq. \eqref{eq4.7} has at least four distinct solutions when considering rotational symmetry. More specifically, solutions III-VI in Fig. \ref{New2figure2} can be rotated at specific angles to transform into one another, a similar characteristic symmetry shared by solutions VII-X.

\begin{figure}[H]
	\begin{center}
		\subfigure[$\phi_0$ ]{ \includegraphics[scale=.2]{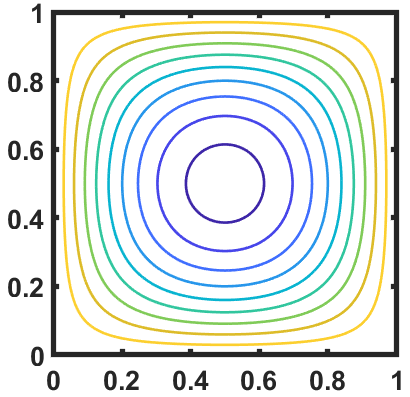}}
		\subfigure[$\phi_1$ ]{\includegraphics[scale=.2]{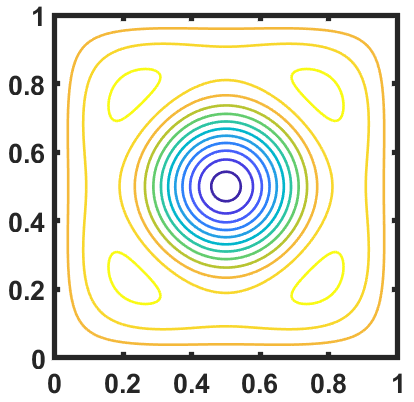}}
		\subfigure[$\phi_2$ ]{ \includegraphics[scale=.2]{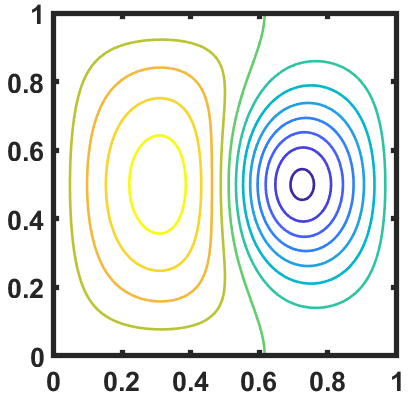}}
		\subfigure[$\phi_3$]{\includegraphics[scale=.2]{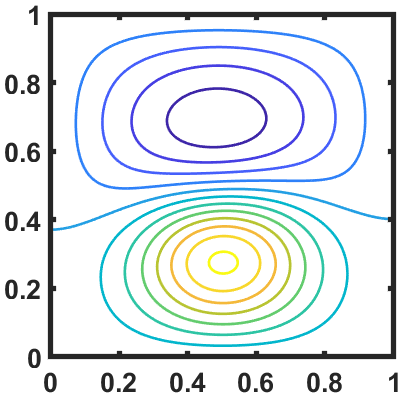}}
		\subfigure[$\phi_4$]{\includegraphics[scale=.2]{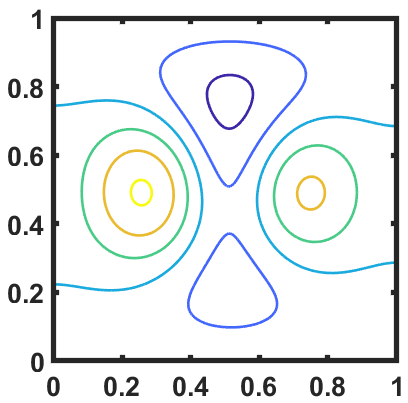}}
	\;\;	\subfigure[$\phi_5$ ]{\includegraphics[scale=.2]{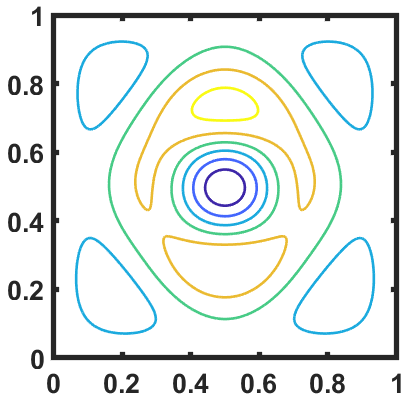}}
		\subfigure[$\phi_6$ ]{\includegraphics[scale=.2]{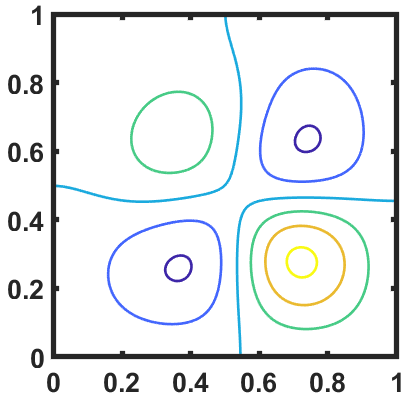}}
		\subfigure[$\phi_7$]{\includegraphics[scale=.2]{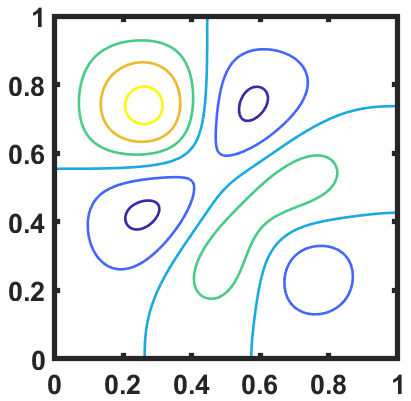}}	
		\subfigure[$\phi_8$]{\includegraphics[scale=.2]{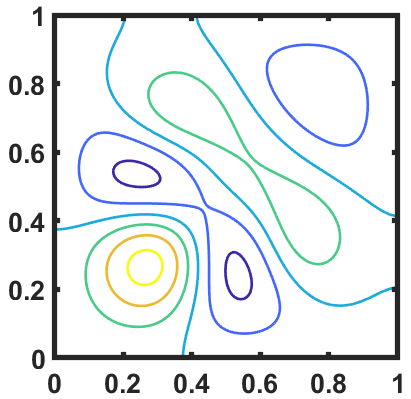}}
		\subfigure[$\phi_9$]{\includegraphics[scale=.2]{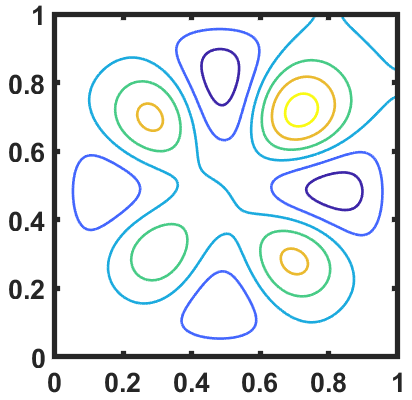}}
		\caption{Basis functions computed by our algorithm for solving \eqref{eq4.7}.}\label{New2figure1}
	\end{center}
\end{figure}

\begin{figure}[H]
	\begin{center}
		\subfigure[I ]{ \includegraphics[scale=.4]{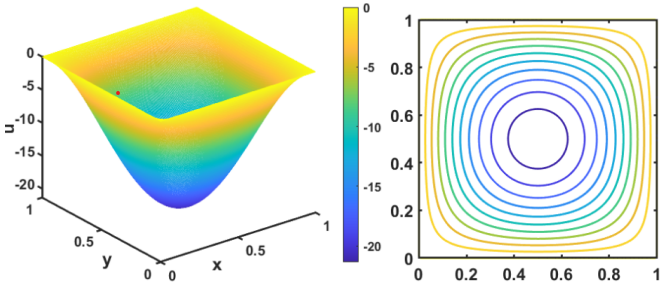}}\;\;
		\subfigure[II ]{\includegraphics[scale=.4]{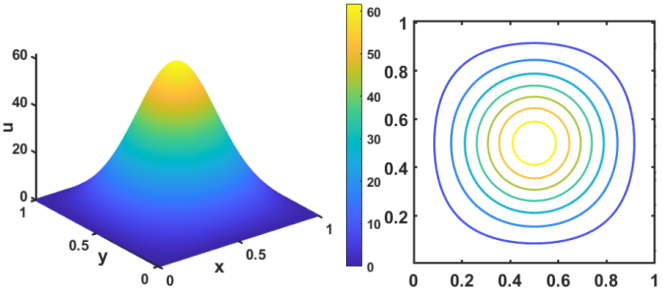}}\\
		
		\subfigure[III ]{ \includegraphics[scale=.4]{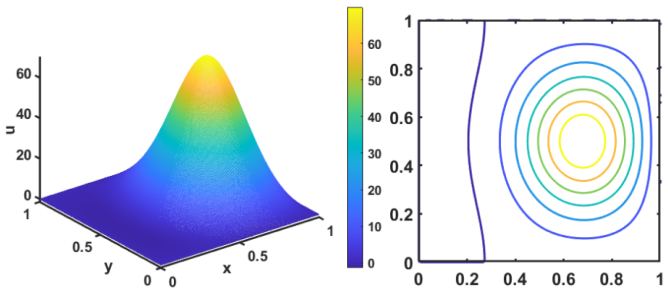}}\;\;
		\subfigure[IV ]{\includegraphics[scale=.4]{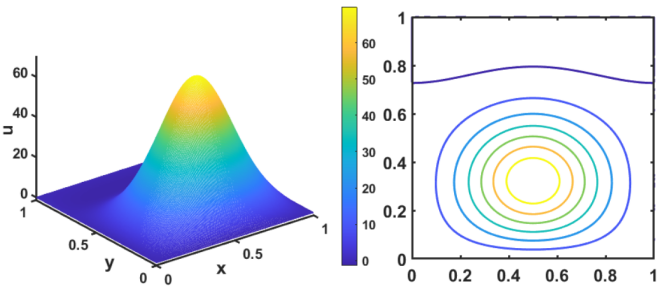}}\\
		
		\subfigure[V ]{\includegraphics[scale=.4]{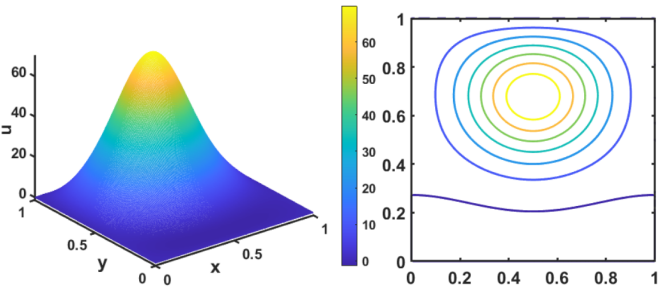}}\;\;
		\subfigure[VI ]{\includegraphics[scale=.4]{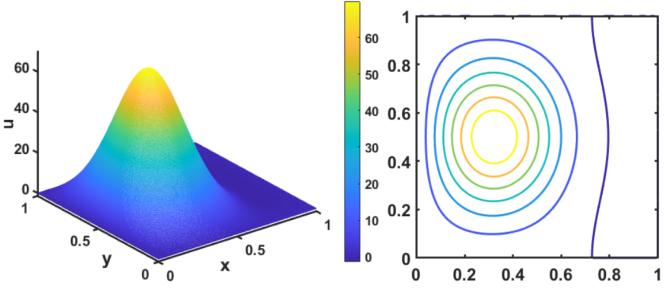}}\\

		\subfigure[VII ]{ \includegraphics[scale=.4]{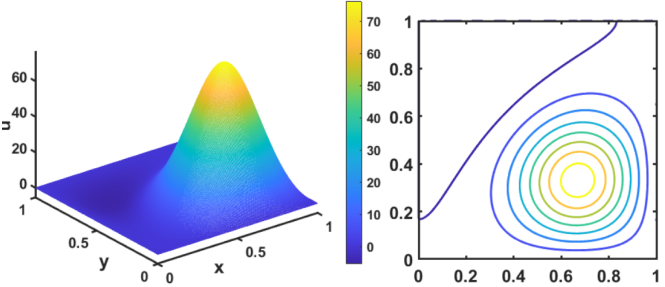}}\;\;
        \subfigure[VIII]{\includegraphics[scale=.4]{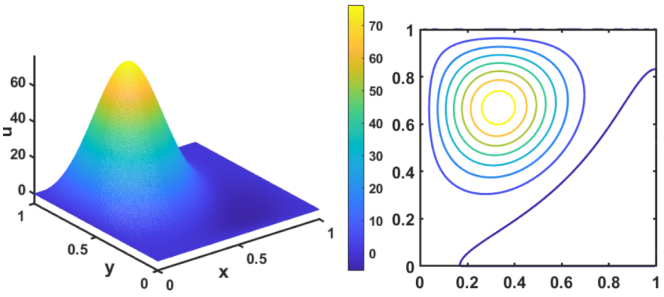}}\\

        \subfigure[IX ]{ \includegraphics[scale=.4]{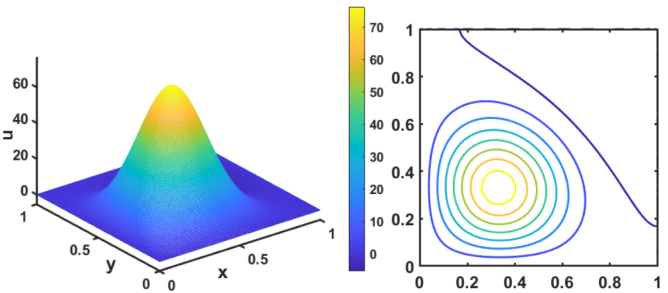}}\;\;
        \subfigure[X ]{\includegraphics[scale=.4]{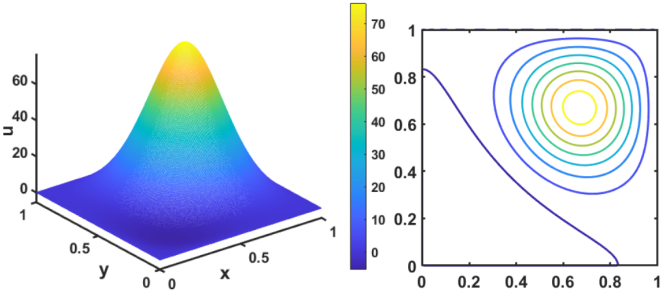}}\\
		  \caption{Multiple solutions of \eqref{eq4.7} by using our algorithm.}\label{New2figure2}
	\end{center}
\end{figure}

The adaptive basis functions generated by our algorithm are shown in Fig. \ref{New2figure1}, and the corresponding coefficients ${\tilde{\alpha}_{9, i}}$ are listed in Table \ref{New2table1}. The computational time required is 157 seconds, with a residual error of less than $10^{-10}$.

\begin{table}[h]
 \renewcommand{\arraystretch}{1.2}
 \setlength\tabcolsep{1.3pt}
	\centering
	\caption{The coefficients associated with various basis functions shown in Fig. \ref{New2figure1}, the residual error and computational time of our algorithm for solving Eq. \eqref{eq4.7}.}\label{New2table1}
	\scalebox{0.7}{
	\begin{tabular}{|c|cccccccccc|c|c|}
		\hline
		\diagbox{Sol.ind}{Coefs}& $\tilde{\alpha}_{9,0}$ & $\tilde{\alpha}_{9,1}$ & $\tilde{\alpha}_{9,2}$ & $\tilde{\alpha}_{9,3}$ & $\tilde{\alpha}_{9,4}$ & $\tilde{\alpha}_{9,5}$ & $\tilde{\alpha}_{9,6}$& $\tilde{\alpha}_{9,7}$& $\tilde{\alpha}_{9,8}$& $\tilde{\alpha}_{9,9}$ &Residual & Time(s) \\ \hline
		
		I & 22.49 & {\color{red}5.36e-15} & {\color{red}-2.33e-15} & {\color{red}2.04e-15} & {\color{red}-3.45e-15} & {\color{red}1.34e-15} & {\color{red}-4.81e-16} & {\color{red}2.27e-15}  & {\color{red}2.57e-15}& {\color{red}-2.12e-16}&1.30e-11&\multirow{10}{*}{157.46}\\
		II &  -44.14 & -11.71 &  {\color{red}2.31e-13} & {\color{red}2.56e-13} & {\color{red}1.07e-13} & {\color{red}-3.18e-14} & {\color{red}-5.92e-15} &{\color{red}5.72e-15}  &{\color{red}-4.77e-15}&{\color{red}-2.22e-15}&3.49e-11&\\
		III &  -36.72& -5.42 &  -27.06 & {\color{red}7.65e-13} &  {\color{red}1.61e-13} & {\color{red}-9.03e-14} &   {\color{red}4.57e-13} & {\color{red}8.59e-14} &{\color{red}6.97e-14}&{\color{red}-6.86e-14} & 4.87e-11&\\
		IV &-36.72& -5.42  & 1.03 & 27.05 & {\color{red}8.20e-15}  & {\color{red}-1.99e-14} & {\color{red}3.81e-14} & {\color{red}8.42e-15} &{\color{red}-2.29e-14}&{\color{red}7.54e-15} &3.22e-11&\\
		V &-36.72& -5.42 &  24.76 & -1.97 & 10.76 &  {\color{red}-4.08e-14} & {\color{red}-1.75e-13} & {\color{red}-1.40e-13}  &{\color{red}4.44e-14}&{\color{red}4.61e-15}&3.72e-11&\\
		VI & -36.72& -5.42  & 1.03&- 24.81& -9.47& 5.11 & {\color{red}2.74e-14}  & {\color{red}2.15e-14}  &{\color{red}6.96e-15}&{\color{red}-6.24e-15}& 3.75e-11&\\
		VII & -30.61 & -0.69 &  -22.30 &23.16 & -0.25 & -0.98 & 12.70 & {\color{red}9.68e-16}  &{\color{red}5.30e-15} &{\color{red}-2.30e-15} &5.41e-11&\\
		VIII &-30.61 & -0.69 &22.18&-23.03&0.86&3.40&8.27&9.63&{\color{red}3.65e-15}  &{\color{red}-1.59e-15} & 5.40e-11&\\
		IX&-30.61 & -0.69 &22.18&21.47&8.99&-0.98&-10.31&-4.74&5.70&{\color{red}2.21e-15}&3.79e-11& \\
		X&-30.61 & -0.69 &-22.30&-21.35&-8.38&3.40&-10.31&-4.74&-4.16&3.90&3.77e-11& \\
		\hline
	\end{tabular}
}
\end{table}

\begin{figure}[H]
	\begin{center}
		\subfigure[I ]{ \includegraphics[scale=.25]{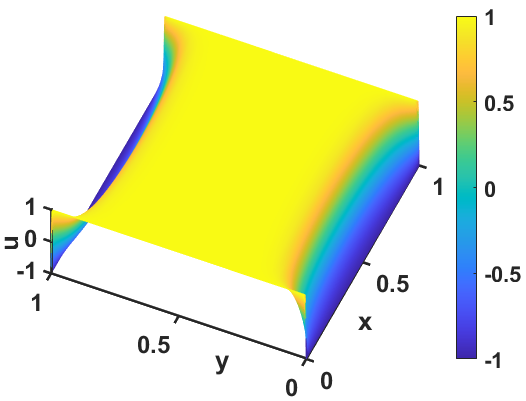}}\;
		\subfigure[II ]{\includegraphics[scale=.25]{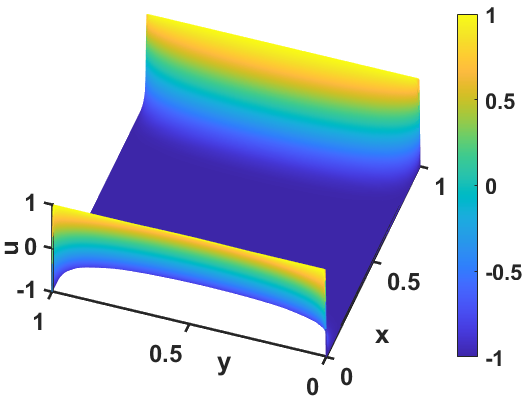}}\;
		\subfigure[III ]{ \includegraphics[scale=.25]{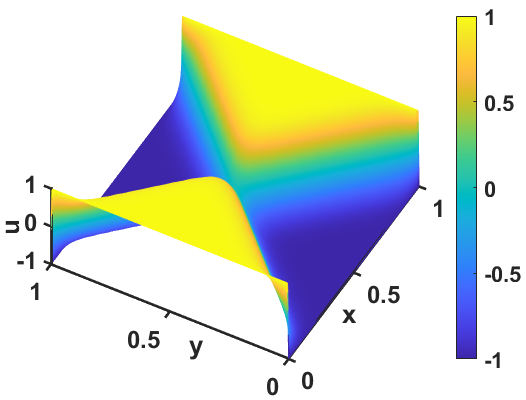}}
		  \caption{Multiple solutions of \eqref{eq4.8} by using our algorithm with $\epsilon=1.6\times 10^{-3}$.}\label{New2figure3}
	\end{center}
\end{figure}

\begin{figure}[H]
	\begin{center}
		\subfigure[$\phi_0$ ]{ \includegraphics[scale=.25]{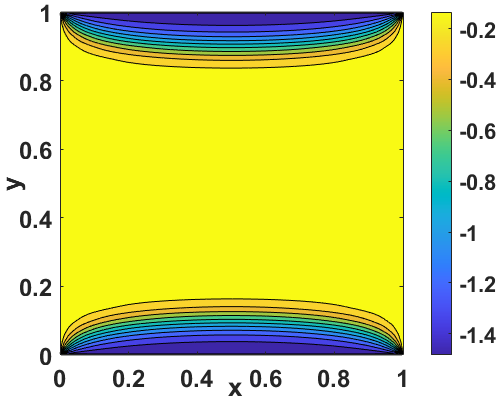}}\;
		\subfigure[$\phi_1$ ]{\includegraphics[scale=.25]{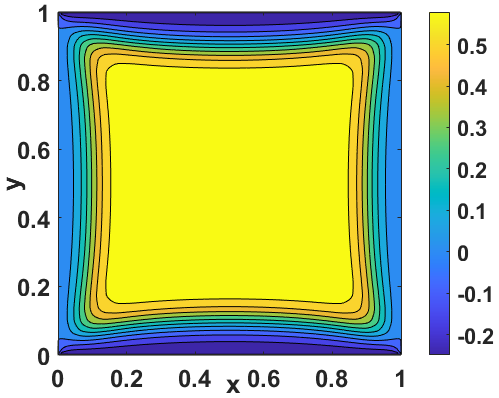}}\;
		\subfigure[$\phi_2$ ]{ \includegraphics[scale=.25]{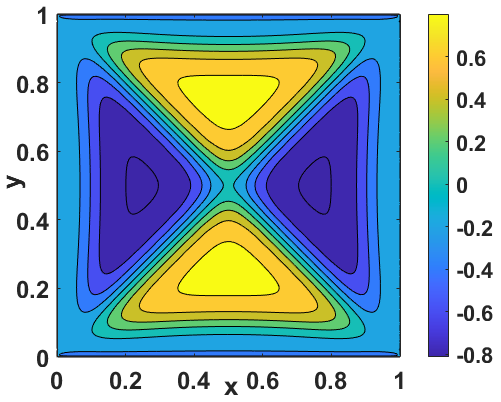}}
	
		  \caption{Basis functions computed by our algorithm for solving \eqref{eq4.8} with $\epsilon=1.6\times 10^{-3}$.}\label{New2figure4}
	\end{center}
\end{figure}

\begin{figure}[H]
	\begin{center}
		\subfigure[I ]{ \includegraphics[scale=.2]{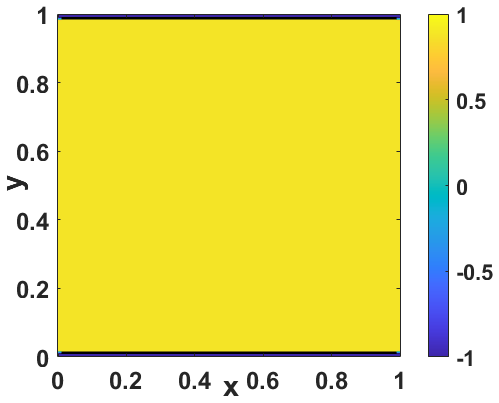}}\;
		\subfigure[II ]{\includegraphics[scale=.2]{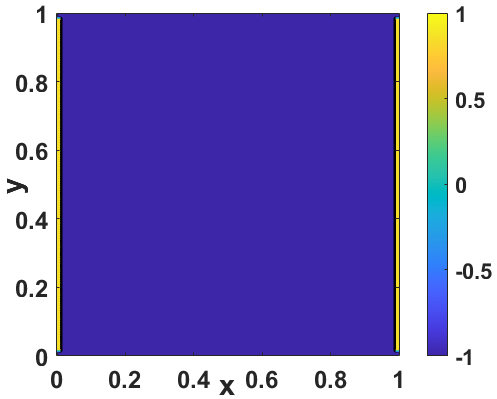}}\;
		\subfigure[III ]{ \includegraphics[scale=.2]{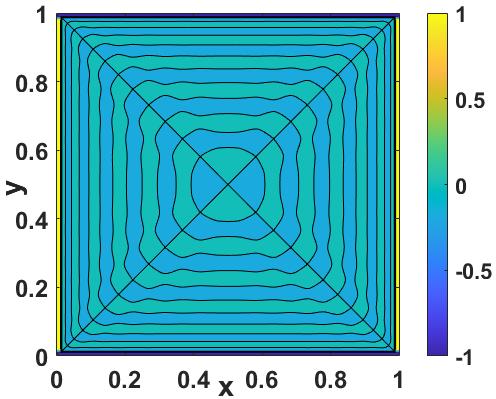}}\;
		\subfigure[IV ]{\includegraphics[scale=.2]{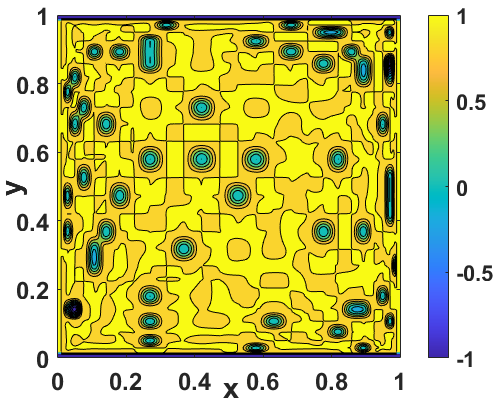}}\\
		\;\subfigure[V ]{\includegraphics[scale=.2]{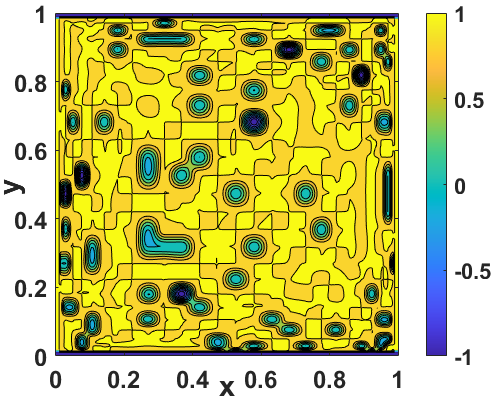}}\;
		\subfigure[VI ]{\includegraphics[scale=.2]{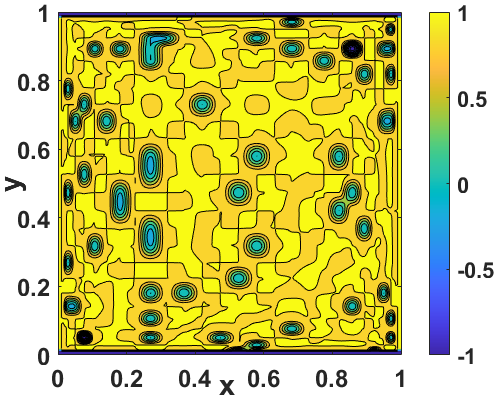}}\;
		\subfigure[VII ]{ \includegraphics[scale=.2]{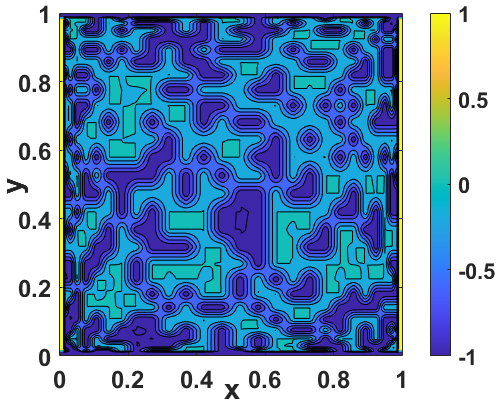}}\;
        \subfigure[VIII]{\includegraphics[scale=.2]{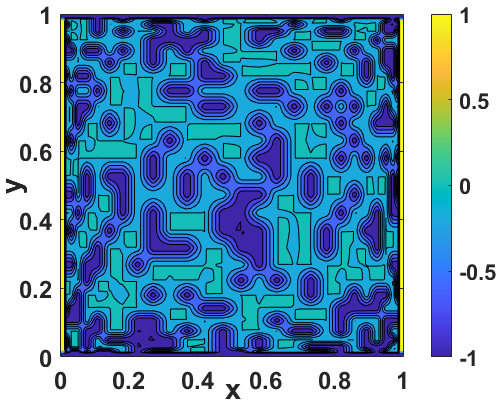}}\\
       \subfigure[IX ]{ \includegraphics[scale=.2]{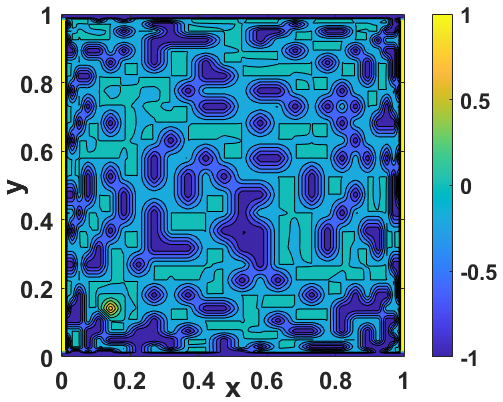}}\;
        \subfigure[X ]{\includegraphics[scale=.2]{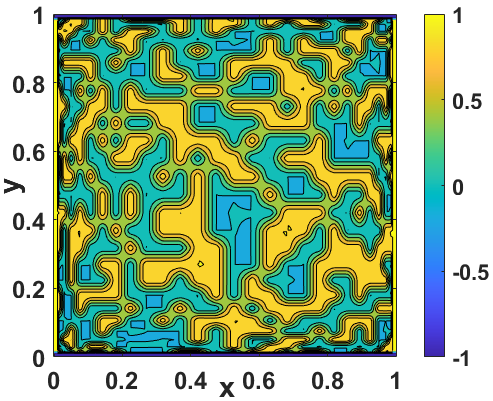}}\;
        \subfigure[XI ]{ \includegraphics[scale=.2]{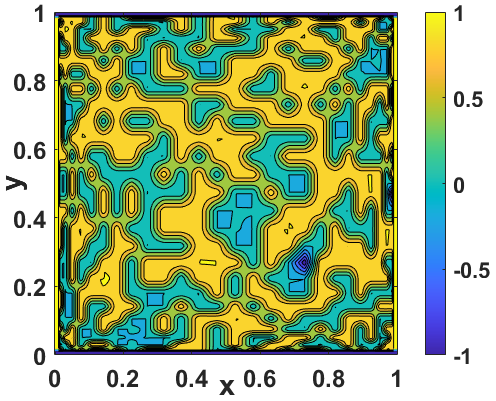}}\;
        \subfigure[XII ]{\includegraphics[scale=.2]{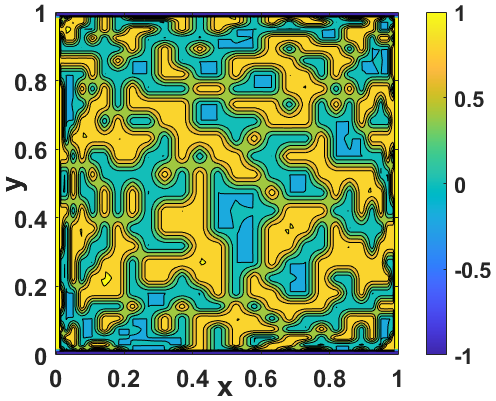}}\\
		  \caption{twelve solutions of \eqref{eq4.8} with $\epsilon=1.0\times 10^{-6}$.}\label{New2figure5}
	\end{center}
\end{figure}

\begin{figure}[H]
	\begin{center}
		\subfigure[I ]{ \includegraphics[scale=.3]{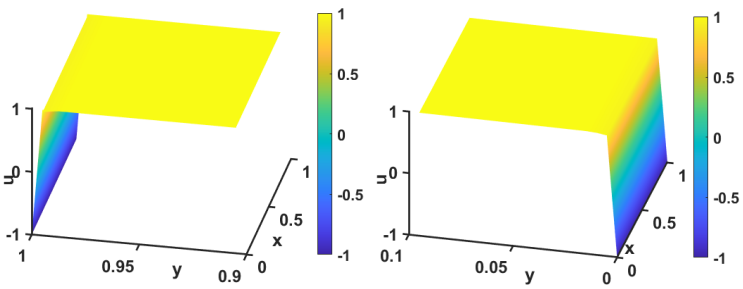}}\;
		\subfigure[III ]{\includegraphics[scale=.35]{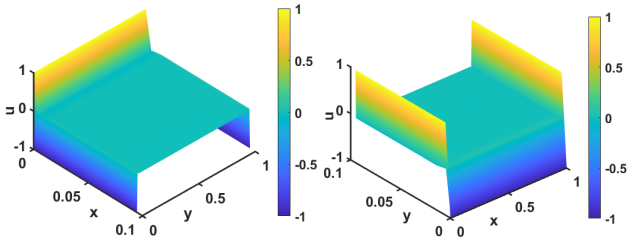}}\;
	\end{center}
 \caption{The behavior of solutions I and III near the boundaries with $\epsilon=1.0\times10^{-6}$.}\label{New2figure6}
\end{figure}

\begin{figure}[H]
	\begin{center}
		\subfigure[$\phi_0$ ]{ \includegraphics[scale=.2]{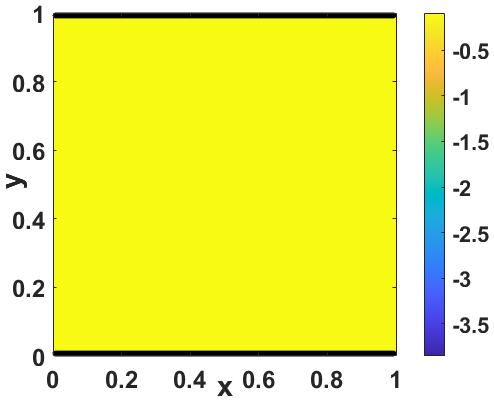}}
		\subfigure[$\phi_1$ ]{\includegraphics[scale=.2]{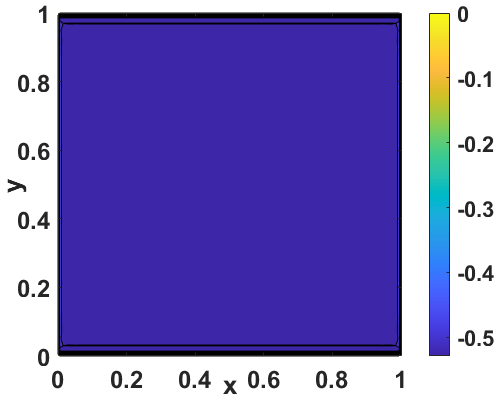}}
		\subfigure[$\phi_2$]{ \includegraphics[scale=.2]{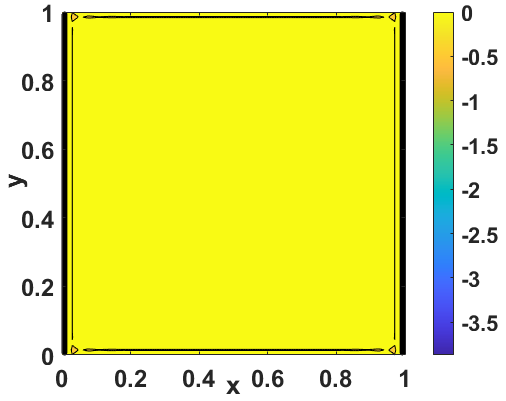}}
		\subfigure[$\phi_3$]{\includegraphics[scale=.2]{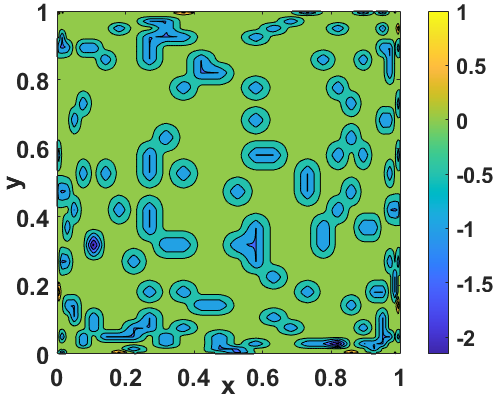}}\\
	\;	\subfigure[$\phi_4$]{\includegraphics[scale=.2]{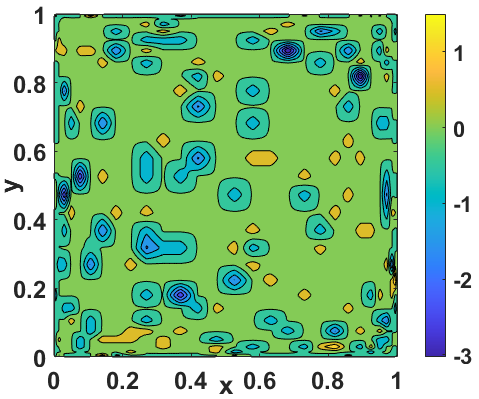}}\;\;
		\subfigure[$\phi_5$]{\includegraphics[scale=.2]{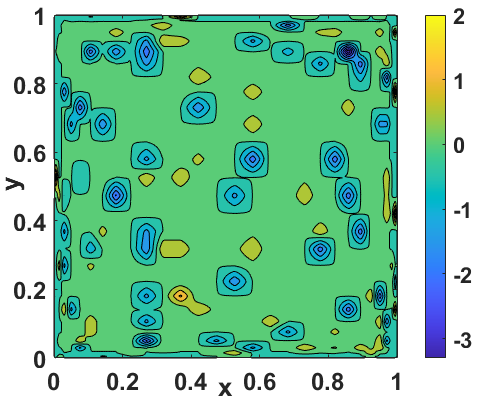}}\;\;
		\subfigure[$\phi_6$]{ \includegraphics[scale=.2]{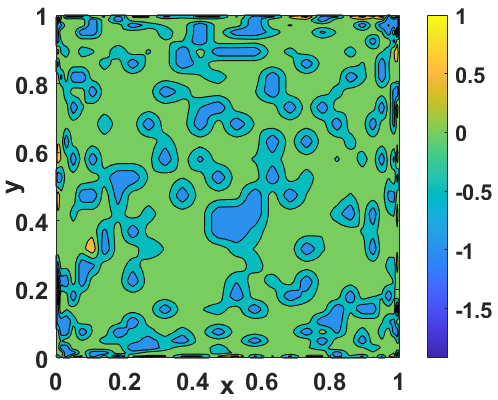}}
        \subfigure[$\phi_7$]{\includegraphics[scale=.2]{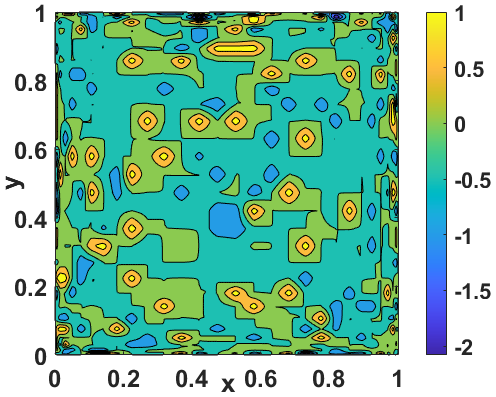}}\\
        \subfigure[$\phi_8$ ]{ \includegraphics[scale=.2]{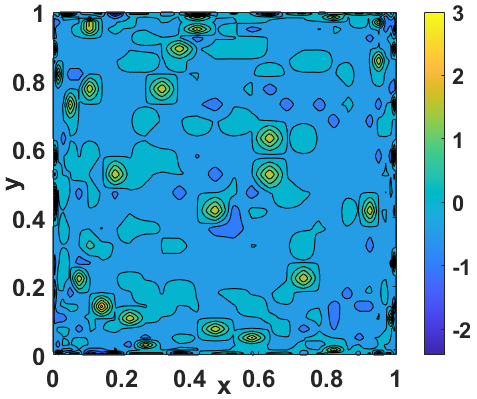}}\;\;
        \subfigure[$\phi_9$]{\includegraphics[scale=.2]{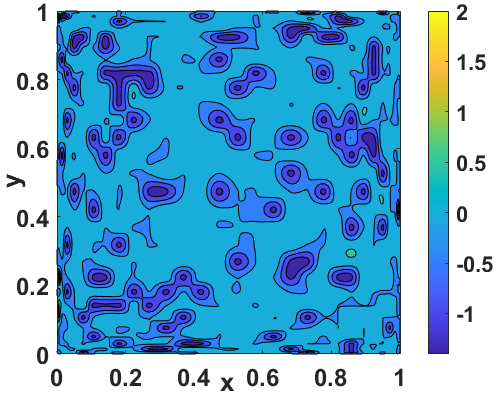}}\;
        \subfigure[$\phi_{10}$ ]{ \includegraphics[scale=.2]{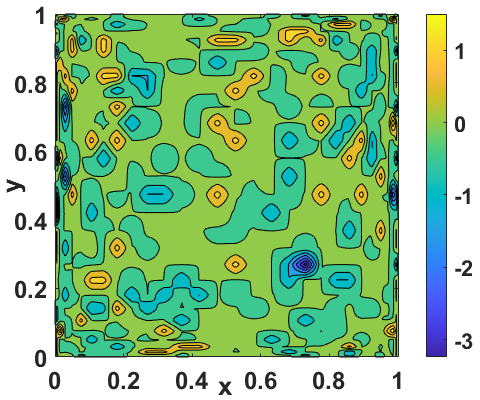}}\;
        \subfigure[$\phi_{11}$ ]{\includegraphics[scale=.2]{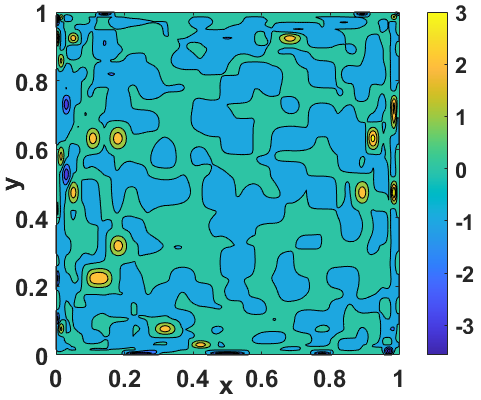}}\\
		  \caption{Basis functions of \eqref{eq4.8} with $\epsilon=1.0\times 10^{-6}$.}\label{New2figure7}
	\end{center}
\end{figure}

\begin{figure}[H]
	\begin{center}
		\subfigure[I ]{ \includegraphics[scale=.2]{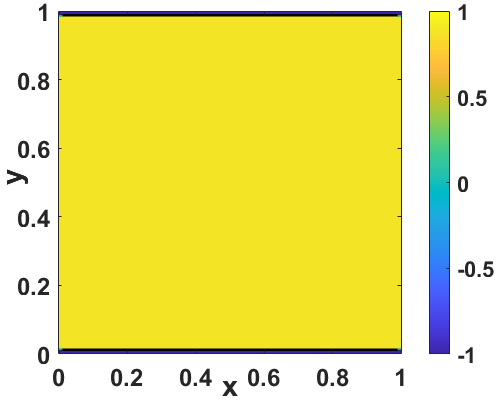}}\;
		\subfigure[II ]{\includegraphics[scale=.2]{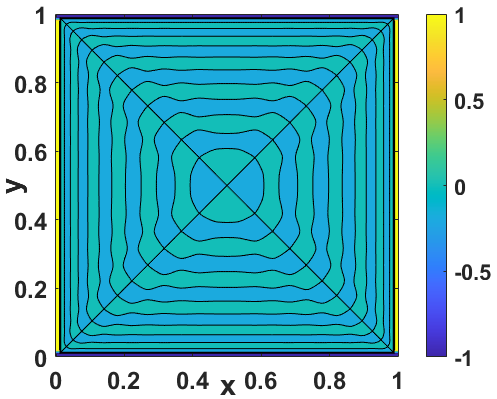}}\;	
		\subfigure[III ]{ \includegraphics[scale=.2]{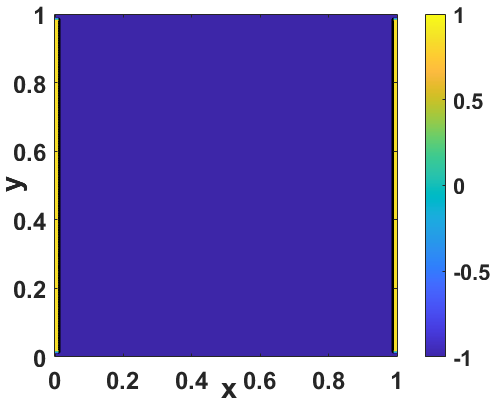}}\;
		\subfigure[IV ]{\includegraphics[scale=.2]{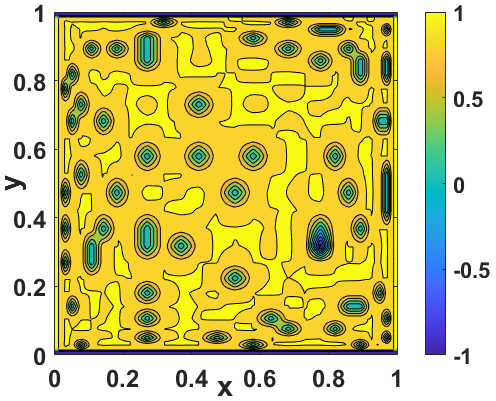}}\\
	\;	\subfigure[V ]{\includegraphics[scale=.2]{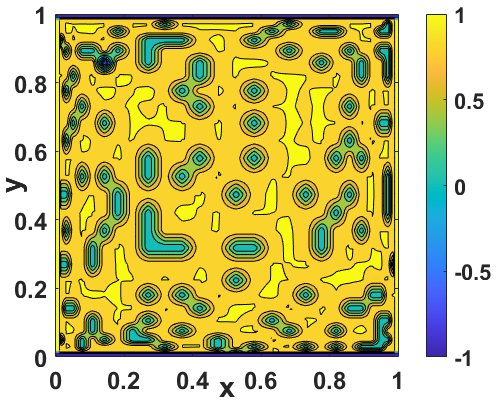}}\;
		\subfigure[VI ]{\includegraphics[scale=.2]{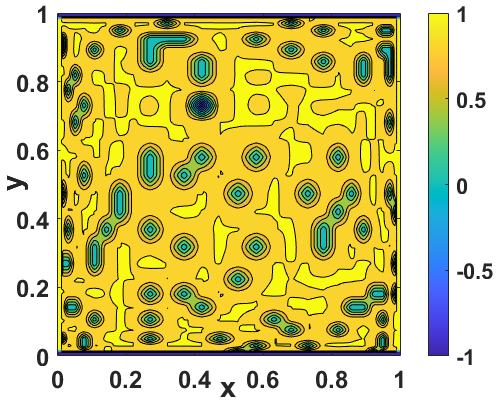}}\;
        \subfigure[VII ]{\includegraphics[scale=.2]{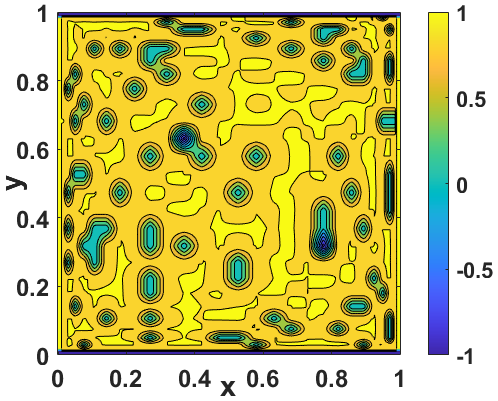}}\;
       \subfigure[VIII ]{ \includegraphics[scale=.2]{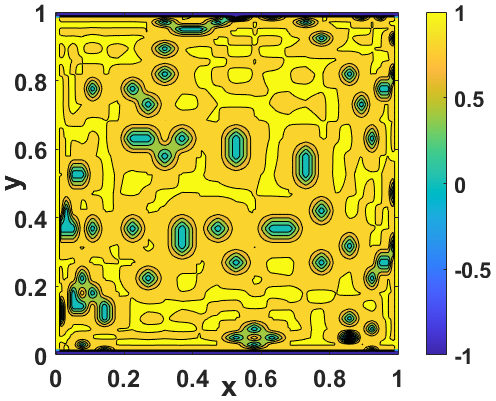}}\\	
		\subfigure[IX ]{ \includegraphics[scale=.2]{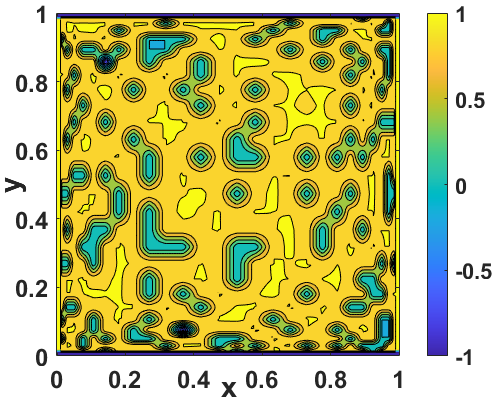}}\;
		\subfigure[X ]{\includegraphics[scale=.2]{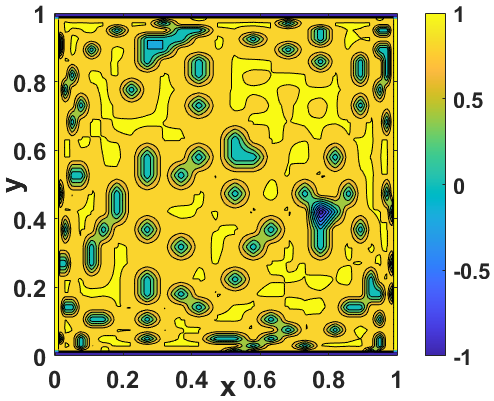}}\;
		\subfigure[XI ]{ \includegraphics[scale=.2]{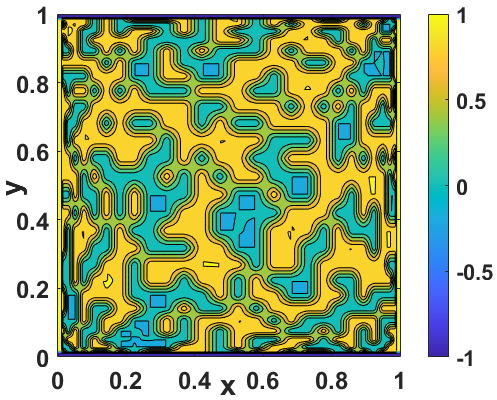}}\;
        \subfigure[XII]{\includegraphics[scale=.2]{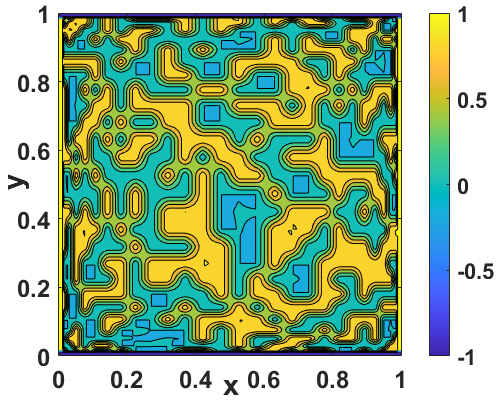}}\\
                \subfigure[XIII ]{ \includegraphics[scale=.2]{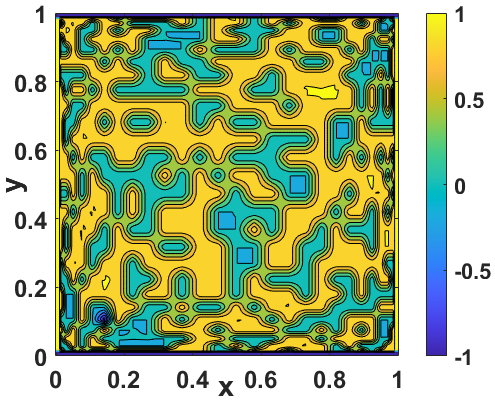}}\;
        \subfigure[XIV ]{\includegraphics[scale=.2]{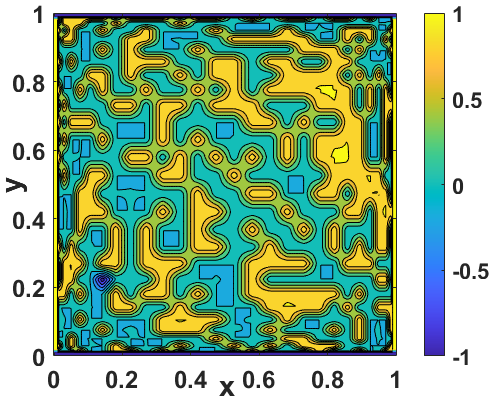}}\;
        \subfigure[XV ]{ \includegraphics[scale=.2]{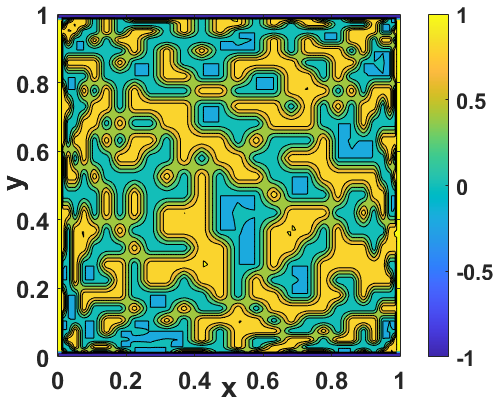}}\;
          \subfigure[XVI ]{\includegraphics[scale=.2]{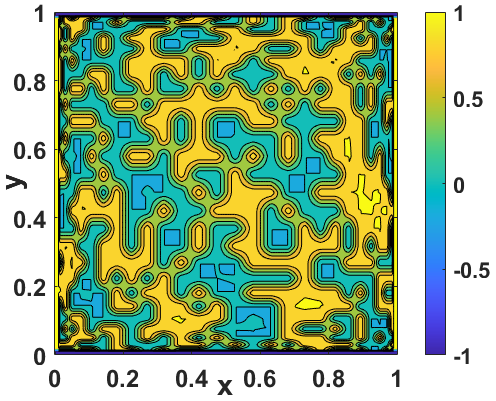}}\\
  \;      \subfigure[XVII ]{\includegraphics[scale=.2]{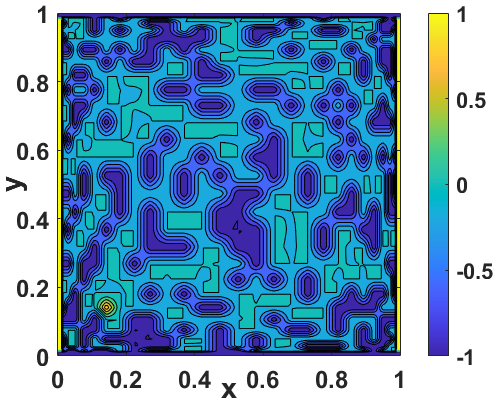}}\;
        \subfigure[XVIII ]{ \includegraphics[scale=.2]{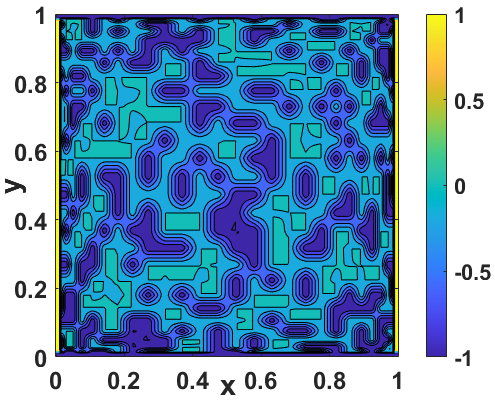}}\;
        \subfigure[XIX ]{\includegraphics[scale=.2]{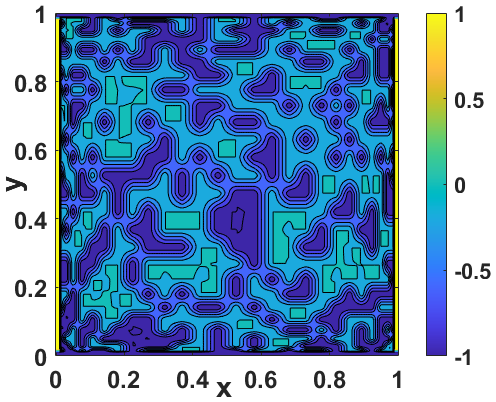}}\;
        \subfigure[XX ]{ \includegraphics[scale=.2]{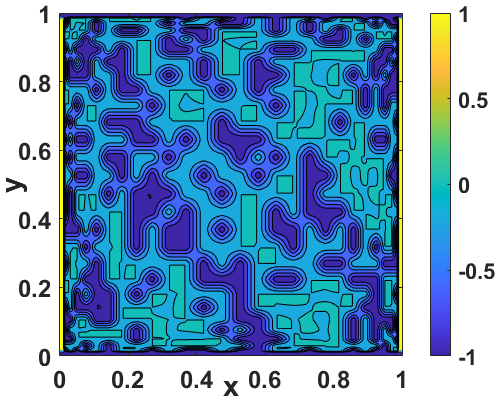}}\\

		  \caption{twenty solutions of \eqref{eq4.8} with $\epsilon=1.0\times 10^{-8}$.}\label{New2figure8}
	\end{center}
\end{figure}

\emph{Example 5.} We consider the steady-state Allen-Cahn equation\cite{2022Two}, described by the following equations:
\begin{equation}\label{eq4.8}
\begin{cases}\vspace{0.15cm}
-\epsilon\Delta{u}+u^3-u=0, \;\;\;\;\;\;\;(x,y) \in \Omega=(0,1)^2,\\
u(0,y)=u(1,y)=1,\;\;\;\;\;\;\;\;y\in(0,1),\\
u(x,0)=u(x,1)=-1,\;\;\;\;\;x\in(0,1).
\end{cases}
\end{equation}
Here, $\epsilon$ serves as a parameter that characterizes the balance between free surface tension and the potential term in the free energy. Here our main concern is that multiple solutions are computed when $\epsilon \rightarrow 0$. To better show our numerical results, three cases (i.e. $\epsilon = 1.6\times 10^{-3}, 1.0\times 10^{-6}$ and $1.0\times 10^{-8}$) are mainly considered. When $\epsilon = 1.6\times 10^{-3}$, three solutions for \eqref{eq4.8} are computed, and plotted in Fig. \ref{New2figure3}. The adaptive basis functions are depicted in Fig. \ref{New2figure4} and corresponding coefficients, denoted as ${\tilde{\alpha}_{2, i}}$,  are documented in Table \ref{New2table2}. Next, we make a test with $\epsilon=1.0\times 10^{-6}$, many solutions can be computed. For simplicity, here twelve solutions are shown in Fig.\ref{New2figure5}. Their adaptive basis functions are presented in Fig.\ref{New2figure7}, and the associated coefficients, denoted as ${\tilde{\alpha}_{11, i}}$, are listed in Table \ref{New2table3}. To check the behavior of these solutions near the boundaries, the pictures of solutions I and III near the boundaries are presented in Fig.\ref{New2figure6}, where the boundary layer can be observed. Finally, when $\epsilon=1.0\times 10^{-8}$, more numerical solutions can be found compared with the case $\epsilon=1.0\times 10^{-6}$. For simplicity, here twenty solutions are shown in Fig.\ref{New2figure8}. In addition, in some published literature, e.g. \cite{2022XinHao}, it has been pointed out that when $\epsilon \rightarrow 0$, there are more and more solutions to \eqref{eq4.8}, which is consistent with our numerical results.

\begin{table}[h]
	\renewcommand{\arraystretch}{1.2}
	\setlength\tabcolsep{1.3pt}
	\centering
	\caption{Numerical results of \eqref{eq4.8} with $\epsilon=1.6\times 10^{-3}$. }\label{New2table2}
	\begin{tabular}{|c|ccc|c|c|}
		\hline
		\diagbox{Sol.ind}{Coefs} & $\tilde{\alpha}_{2,0}$ & $\tilde{\alpha}_{2,1}$ & $\tilde{\alpha}_{2,2}$  &Residual & Time(s) \\ \hline
		I & 1.3358 & {\color{red}-6.1862e-17} & {\color{red}-3.4779e-17}  & 5.2625e-14 &\multirow{3}{*}{32.74}\\
		II &  1.8557 & -3.0235 &  {\color{red}2.0258e-16} &3.9191e-14  &~\\
		III &  1.8020 & -1.4763 &  -1.1867 & 6.7391e-14 &~ \\

		\hline
	\end{tabular}
\end{table}


\begin{table}[h]
	\renewcommand{\arraystretch}{1.2}
	\setlength\tabcolsep{1.3pt}
	\centering
	\caption{Numerical results of \eqref{eq4.8}  with $\epsilon=1.6\times 10^{-6}$. }\label{New2table3}
	\scalebox{0.6}{
	\begin{tabular}{|c|cccccccccccc|c|c|}
		\hline
		\diagbox{Sol.ind}{Coefs}& $\tilde{\alpha}_{11,0}$ & $\tilde{\alpha}_{11,1}$ & $\tilde{\alpha}_{11,2}$ & $\tilde{\alpha}_{11,3}$ & $\tilde{\alpha}_{11,4}$ & $\tilde{\alpha}_{11,5}$ & $\tilde{\alpha}_{11,6}$& $\tilde{\alpha}_{11,7}$& $\tilde{\alpha}_{11,8}$& $\tilde{\alpha}_{11,9}$& $\tilde{\alpha}_{11,11}$& $\tilde{\alpha}_{11,11}$ &Residual & Time(s) \\ \hline
		
		I & 0.01 & {\color{red}-1.54e-16} & {\color{red}-3.39e-17} & {\color{red}-1.44e-17} & {\color{red}1.71e-17} & {\color{red}-8.18e-19} & {\color{red}1.58e-17} & {\color{red}-1.37e-17}  & {\color{red}7.23e-18}& {\color{red}2.93e-17}&{\color{red}-9.60e-18}&{\color{red}1.17e-17}&8.88e-16&\multirow{12}{*}{67.28}\\
		II &  0.44 & 3.97 &  {\color{red}-1.79e-18} & {\color{red}2.73e-17} & {\color{red}1.23e-17} & {\color{red}3.26e-17} & {\color{red}-2.29e-17} &{\color{red}1.48e-17}  &{\color{red}3.43e-17}&{\color{red}2.70e-17}&{\color{red}5.42e-17}&{\color{red}-1.04e-17}&4.44e-16&\\
		III &  0.21& 1.98 &  0.02 & {\color{red}-3.50e-17} &  {\color{red}5.84e-18} & {\color{red}2.36e-17} & {\color{red}-1.94e-17} & {\color{red}3.90e-18} &{\color{red}1.66e-17}&{\color{red}1.45e-17} &{\color{red}2.17e-17} &{\color{red}1.42e-17} & 2.77e-16&\\
		IV &0.24& 0.40  & 0.24 & 0.80 & {\color{red}-3.20e-17}  & {\color{red}6.08e-17} & {\color{red}5.05e-17} & {\color{red}9.21e-19} &{\color{red}8.96e-17}&{\color{red}1.97e-17} &{\color{red}-2.95e-17} &{\color{red}4.62e-18} &6.66e-16&\\
		V &0.23& 0.27 &  0.22 & 0.35 & 0.60 &  {\color{red}4.07e-16} & {\color{red}3.17e-17} & {\color{red}5.11e-17}  &{\color{red}2.23e-17}&{\color{red}-3.19e-16}&{\color{red}1.29e-16}&{\color{red}5.75e-16}&6.66e-16&\\
		VI & 0.23& 0.18  & 0.21&0.21& 0.19& 0.46 & {\color{red}1.04e-16}  & {\color{red}1.49e-17}  &{\color{red}-8.85e-17}&{\color{red}-5.45e-16}&{\color{red}2.16e-16}&{\color{red}8.48e-16}& 6.66e-16&\\
		VII & 0.31 & 2.72 &  -0.23 &0.55 & 0.17 & 0.07 & 0.78 & {\color{red}2.02e-17}  &{\color{red}-3.11e-17} &{\color{red}-2.07e-19} &{\color{red}-1.15e-17} &{\color{red}-4.24e-18} &4.44e-16&\\
		VIII &0.26 & 2.57 &-0.17&0.63&0.18&0.07&0.46&0.45&{\color{red}6.27e-17}  &{\color{red}-4.82e-17} &{\color{red}-5.38e-18} &{\color{red}1.05e-17} & 4.44e-16&\\
		IX&0.31 & 2.52 &0.03&0.65&0.19&0.08&0.36&0.30&0.34&{\color{red}1.28e-17}&{\color{red}-9.10e-17}&{\color{red}9.60e-18}&4.44e-16& \\
		X&0.12 & 1.07 &-0.05&0.41&0.12&0.04&0.59&0.02&-0.02&0.67&{\color{red}1.80e-17}&{\color{red}5.61e-17}&6.66e-16& \\
          XI&0.12 & 0.94 &0.04&0.48&0.14&0.06&0.66&0.03&-0.05&0.36&0.44&{\color{red}4.26e-17}&6.66e-16& \\
          XII&0.16 & 1.05 &0.20&0.43&0.13&0.05&0.59&0.02&-0.03&0.57&0.10&0.28&6.66e-16& \\
		\hline
	\end{tabular}
}
\end{table}

\begin{figure}[H]
	\begin{center}
		\subfigure[I ]{ \includegraphics[scale=.2]{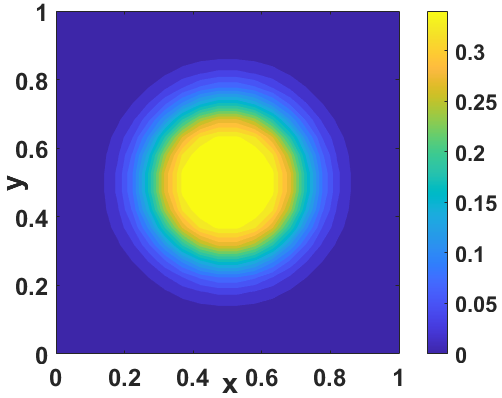}}\;
		\subfigure[II ]{\includegraphics[scale=.2]{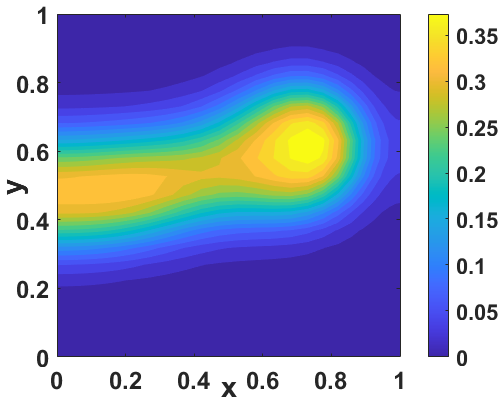}}\;
		\subfigure[III ]{ \includegraphics[scale=.2]{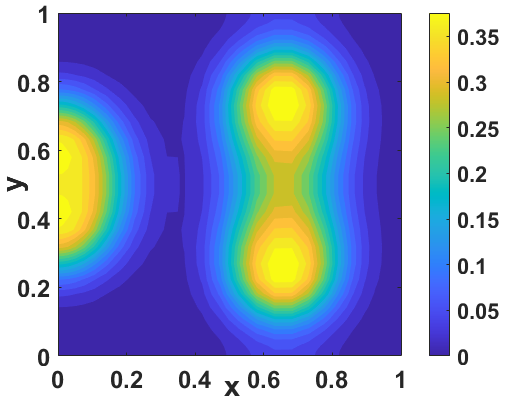}}\;
  	\subfigure[IV ]{ \includegraphics[scale=.2]{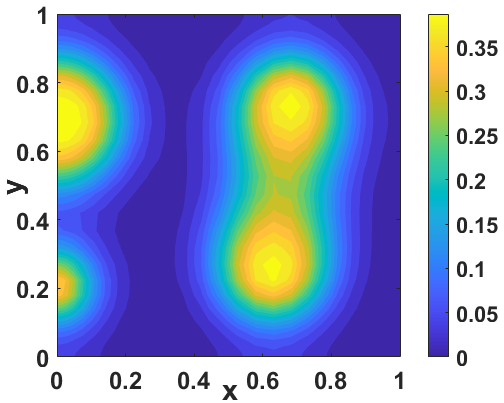}}\;\\
		\subfigure[V ]{\includegraphics[scale=.2]{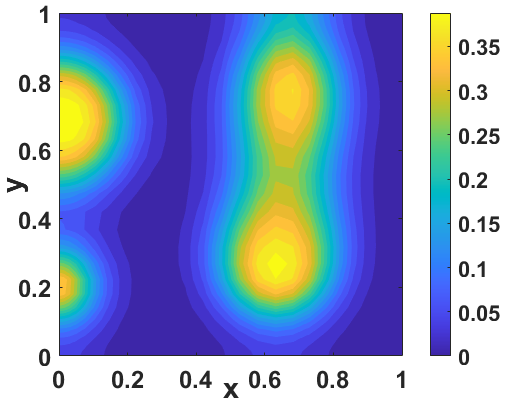}}\;
		\subfigure[VI ]{ \includegraphics[scale=.2]{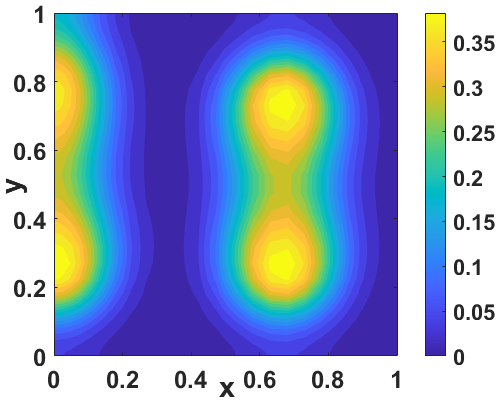}}\;
  	\subfigure[VII ]{ \includegraphics[scale=.2]{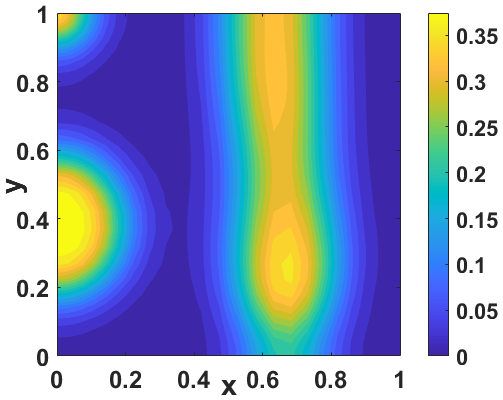}}\;\;
		\subfigure[VIII ]{\includegraphics[scale=.2]{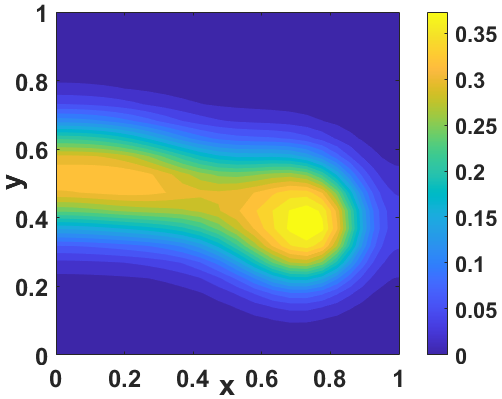}}\;\\
  \caption{Multiple solutions about $A(x,y)$ of \eqref{eq4.9} by using our algorithm.}\label{New3figure1}
	\end{center}
\end{figure}

\begin{figure}[H]
	\begin{center}
		\subfigure[$\phi_0$ ]{ \includegraphics[scale=.2]{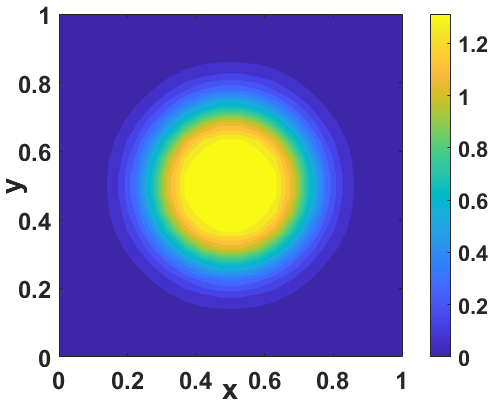}}\;
		\subfigure[$\phi_1$ ]{\includegraphics[scale=.2]{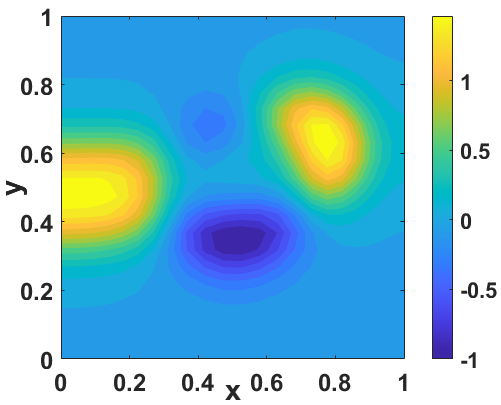}}\;
		\subfigure[$\phi_2$ ]{ \includegraphics[scale=.2]{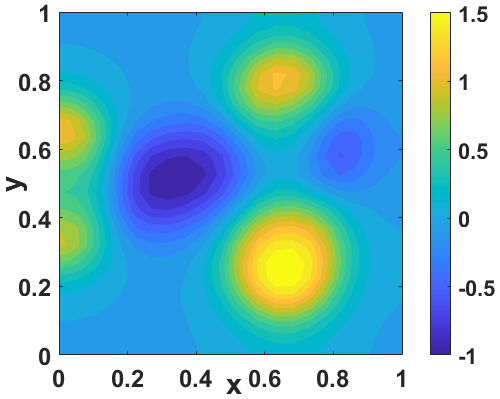}}
  	\subfigure[$\phi_3$ ]{ \includegraphics[scale=.2]{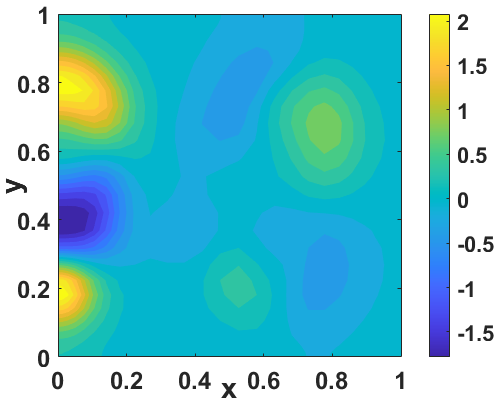}}\;\\
		\subfigure[$\phi_4$]{\includegraphics[scale=.2]{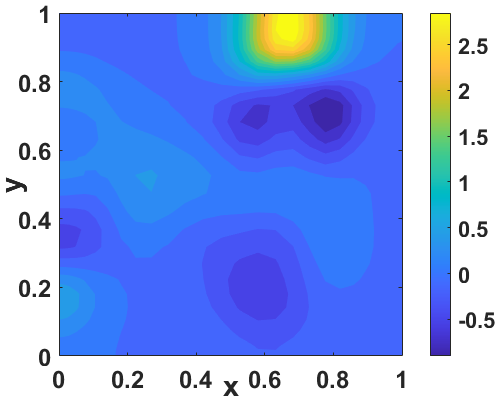}}\;
		\subfigure[$\phi_5$]{ \includegraphics[scale=.2]{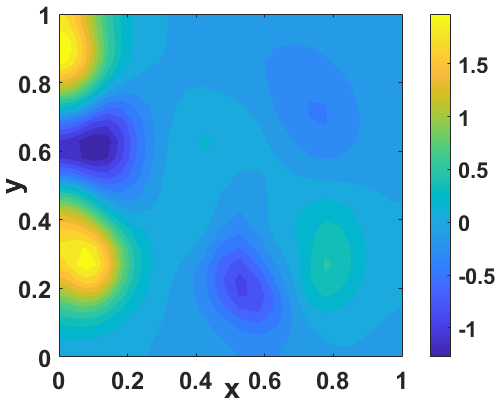}}\;
  	\subfigure[$\phi_6$]{ \includegraphics[scale=.2]{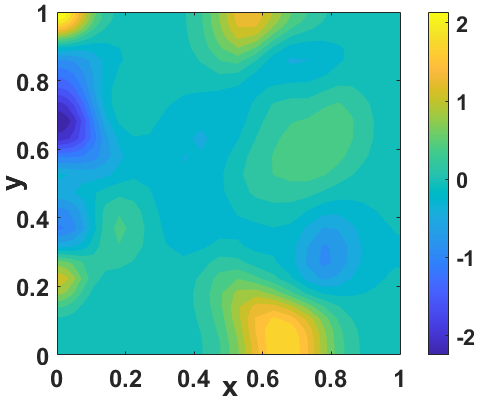}}\;\;
		\subfigure[$\phi_7$ ]{\includegraphics[scale=.2]{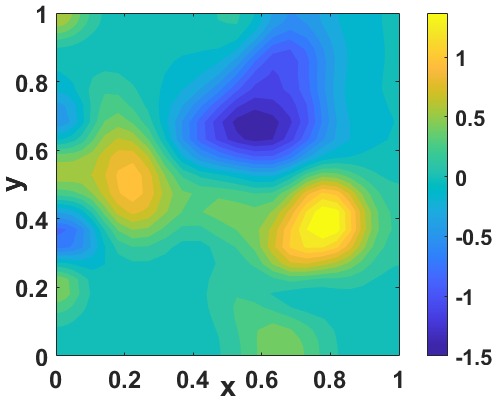}}\;\\
    \caption{Basis functions about $A(x,y)$  for solving \eqref{eq4.9} by our algorithm.}\label{New3figure2}
	\end{center}
\end{figure}

\emph{Example 6.} Our final example explores the steady-state Gray-Scott model \cite{hao2020spatial}, described by the following equations:
\begin{equation}\label{eq4.9}
\begin{cases}\vspace{0.15cm}
D_A\Delta{A}=-SA^2+(\mu+\rho)A, \;\;\;\;\;\;\;{\rm on}\;\; \Omega=(0,1)\times(0,1),\\
D_S\Delta{S}=SA^2-\rho(1-S), \;\;\;\;\;\;\;\;\;\;\;{\rm on}\;\; \Omega=(0,1)\times(0,1),\\
\frac{\partial A}{\partial n}=\frac{\partial S}{\partial n}=0,\;\;\;\;\;\;\;\;\;\;\;\;\;\;\;\;\;\;\;\;\;\;\;\;\;\;\;\;\;\;{\rm on}\;\; \partial\Omega.
\end{cases}
\end{equation}
In this example, we fix the parameters as follows: $D_A = 2.5\times10^{-4}$, $D_S = 5\times10^{-4}$, $\rho = 0.04$, and $\mu = 0.065$. Our algorithm obtains eight solutions about $A(x,y)$, which are graphically presented in Fig. \ref{New3figure1}, alongside their adaptive basis functions in Fig. \ref{New3figure2}. The associated coefficients, denoted as ${\tilde{\alpha}_{7, i}}$, are listed in Table \ref{New3table1}.

\begin{table}[h]
	\renewcommand{\arraystretch}{1.2}
	\setlength\tabcolsep{1.3pt}
	\centering
	\caption{Numerical results about $A(x,y)$ of \eqref{eq4.9} }\label{New3table1}
	\scalebox{0.9}{
	\begin{tabular}{|c|cccccccc|c|}
		\hline
		\diagbox{Sol.ind}{Coefs}& $\tilde{\alpha}_{7,0}$ & $\tilde{\alpha}_{7,1}$ & $\tilde{\alpha}_{7,2}$ & $\tilde{\alpha}_{7,3}$ & $\tilde{\alpha}_{7,4}$ & $\tilde{\alpha}_{7,5}$ & $\tilde{\alpha}_{7,6}$& $\tilde{\alpha}_{7,7}$& Residual  \\ \hline
		
		I & 0.25 & {\color{red}1.23e-14} & {\color{red}-1.32e-14} & {\color{red}-1.01e-15} & {\color{red}4.37e-15} & {\color{red}1.47e-15} & {\color{red}-2.76e-15} & {\color{red}-8.11e-15}  & 2.76e-13\\
		II &  0.21 & 0.21 &  {\color{red}-3.26e-16} & {\color{red}-3.03e-15} & {\color{red}4.53e-15} & {\color{red}2.42e-15} & {\color{red}-3.68e-15} &{\color{red}3.43e-14}  & 2.24e-13\\
		III &  0.17 & 0.16 &  0.24 & {\color{red}1.88e-13} &  {\color{red}-1.83e-14} & {\color{red}-1.43e-13} & {\color{red}-9.06e-14} & {\color{red}1.44e-14} & 2.69e-13\\
		IV &0.17& 0.11  & 0.23 & 0.13 & {\color{red}-3.55e-15}  & {\color{red}5.37e-15} & {\color{red}2.01e-15} & {\color{red}-2.20e-15} &2.32e-13\\
		V &0.16& 0.10 &  0.24 & 0.12 & 0.05 &  {\color{red}-1.78e-14} & {\color{red}-1.69e-14} & {\color{red}2.63e-14} & 2.51e-13\\
		VI & 0.16& 0.15  & 0.24&0.08& 9.16e-04& 0.09 & {\color{red}-1.93e-15}  & {\color{red}1.49e-15} &2.10e-13\\
		VII & 0.16 & 0.14 &  0.23 &-0.02 & 0.07 & 0.08 & 0.09 & {\color{red}-1.45e-15} & 1.93e-13\\
		VIII &0.21 & 0.11 &0.07&-0.03&0.03&4.91e-03&-0.05&0.16&2.52e-13\\
		
		\hline
	\end{tabular}
}
\end{table}

\section{Concluding remarks}\label{sect5}

In this paper, we present an innovative approach for computing multiple solutions of nonlinear differential equations. Our method not only generates multiple initial estimates for solving differential equations with polynomial nonlinearities but also adaptively orthogonal basis functions to solve the discretized nonlinear system. Through a series of numerical experiments, we demonstrate the efficiency and robustness of this newly developed adaptive orthogonal basis method. The convergence analysis of this proposed method will be considered as our future work, although the trust region method can be guaranteed to have quadratic convergence in each iteration (as detailed in the appendix).

$$$$
\textbf{Declarations}
\begin{itemize}
\item[$\bullet$]\textbf{Availability of data and materials:} The data that support the findings of this study are available from the corresponding author upon reasonable request.

\end{itemize}

\section*{Appendix}
We present the detailed process of the trust region method to solve \eqref{eq2.2.4.2}. For this purpose, we introduce a region around the current best solution, and approximate the objective function by a quadratic form which boils down to solving a sequence of trust-region subproblems:
\begin{equation}\label{eq2.2.4.4}
\begin{split}
& \min_{\bs s\in {\mathbb B}_{h_k}} q^{(k)}({\bs s}) := Q(\textbf{\textit{a}}_{k}) + {\bs g}(\textbf{\textit{a}}_{k})^{\top}{\bs s} + \frac{1}{2}{\bs s}^{\top}{\bs G}(\textbf{\textit{a}}_{k}){\bs s},   \quad\quad  k \geq 0,
\end{split}
\end{equation}
where the trust region  ${\mathbb B}_{h_k}:=\{\bs s\in {\mathbb R}^n\,:\,\|\bs s\|\le h_k\}$. When $h_{k}$ is given and ${\bs s}_{k}$ is the minimizer of $q^{(k)}({\bs s})$ in \eqref{eq2.2.4.4}, we can update $\textbf{\textit{a}}_{k+1} = \textbf{\textit{a}}_{k} + {\bs s}_{k}$. Obviously, it is one of the most critical steps to choose a proper $h_k$ at each iteration. Based on a good agreement between $q^{(k)}(s_{k})$ and the objective function value $Q({\bs a}_{k+1})$, we should choose $h_{k}$ as large as possible. To be specific, we define a ratio
\begin{equation}\label{rkval}
r_{k} = \frac{Q(\textbf{\textit{a}}_{k}) - Q(\textbf{\textit{a}}_{k+1})}{q^{(k)}(\bs 0) - q^{(k)}({\bs s}_{k})}.
\end{equation}
The ratio $r_{k}$ is an indicator for expanding and contracting the trust region. If $r_{k}$ is negative, the current value of $Q({\bs a}_{k})$ is less than the new objective value $Q({\bs a}_{k} + {\bs s}_{k})$, consequently the step should be rejected. If $r_{k}$ is close to 1, it means there is a good agreement between the model $q^{(k)}$ and the objective function $Q$ over this step, we can expand the trust region for the next iteration. If $r_{k}$ is close to zero, the trust region should be contracted. Otherwise, we do not alter the trust region at the next iteration. Moreover, the process is also summarized in the following \textbf{algorithm 1}.\vspace{0.4cm}

\begin{center}
\begin{tabular}{l}
\toprule[0.5pt]
\textbf{Algorithm 1} - The trust region method \label{3.15}\\\toprule[0.5pt]
\quad \textbf{Input:} Given $\textbf{\textit{a}}_{0}$, $\epsilon > 0,  0<\delta_1 < \delta_2 < 1, 0 < \tau_1 < 1 < \tau_2$ and $h_{0} = \|\bs g_{0}\|$\\
\quad \textbf{Input:} Initial solution set  $S \leftarrow$ the empty set $\emptyset$\\
\quad \textbf{Output:}  $S$ \\
\, 1:  \quad  \textbf{For} $k = 0,~ 1,~ 2,~ \cdots$\\
\, 2:  \quad\quad Compute ${\bs g}_{k}$ and ${\bs G}_{k}$;\\
\, 3:  \quad\quad If $\|{\bs g}_{k}\| \leq \epsilon$ and $|Q(\textbf{\textit{a}}_{k})| \leq \epsilon$, stop;\\
\, 4:  \quad\quad Approximately solve the subproblem (\ref{eq2.2.4.4}) for ${\bs s}_{k}$;\\
\, 5:  \quad\quad Compute $Q(\textbf{\textit{a}}_{k} + {\bs s}_{k})$ and $r_{k}$. If $r_{k}\geq \delta _1,$ then $\textbf{\textit{a}}_{k+1} = \textbf{\textit{a}}_{k} + {\bs s}_{k}$; Otherwise,\\
\,   \quad\quad\quad set $\textbf{\textit{a}}_{k+1} = \textbf{\textit{a}}_{k}$;\\
\, 7:  \quad\quad If $r_{k} < \delta_1$, then $h_{k+1} = \tau_1 h_{k}$;\\
\,   \quad\quad\quad If $r_{k} > \delta_2$ and $\|{\bs s}_{k}\|$ = $h_{k}$, then $h_{k+1} = \tau_2 h_{k}$;\\
\,   \quad\quad\quad Otherwise, set $h_{k+1} = h_{k}$;\\
\, 8:   \quad \textbf{end}\\
\, 9:   \quad return $S$\\
\toprule[0.5pt]
\end{tabular}
\end{center}
For simplicity, in general we choose $\epsilon = 10^{-13}, \delta_1 = 0.25, \delta_2 = 0.75, \tau_1 = 0.5,$ and $\tau_2 = 2$ throughout the paper. Moreover, in the \textbf{Algorithm 1} (see Line 4), the subproblem \eqref{eq2.2.4.4} needs to be solved. Here the so-called dogleg method (see \cite{sun2006optimization, 0a}) is used to solve it, and the process is as follows: Let $s := {\bs a}_{k} - d_{k}{\bs g}_{k}$, and substituting it into \eqref{eq2.2.4.4} yields
\begin{equation*}
q^{(k)}({\bs a}_{k} - d_{k} {\bs g}_{k}) = Q({\bs a}_{k}) - d_{k}\|{\bs g}_{k}\|^{2}_2 + \frac{1}{2}d^2_{k} {\bs g}^{\top}_{k}{\bs G}_{k}{\bs g}_{k}.
\end{equation*}
Based on the exact line search, the step size $d_{k}$ becomes
\begin{equation*}
d_{k} = \frac{\|{\bs g}_{k}\|^{2}_2}{{\bs g}^{\top}_{k}{\bs G}_{k}{\bs g}_{k}}.
\end{equation*}
Consequently the corresponding step along the steepest descent direction is
\begin{equation*}
{\bs s}^C_{k} = - d_{k}{\bs g}_{k} = -\frac{{\bs g}^{\top}_{k}{\bs g}_{k}}{{\bs g}^{\top}_{k}{\bs G}_{k}{\bs g}_{k}}{\bs g}_{k}.
\end{equation*}
On the other hand, the Newtonian step is
\begin{equation*}
{\bs s}^{N}_{k} = -{\bs G}^{-1}_{k}{\bs g}_{k}.
\end{equation*}
If $\|{\bs s}^{C}_{k}\|_{2} = \|d_{k}{\bs g}_{k}\|_{2} \geq h_{k}$, the solution of \eqref{eq2.2.4.4} can be obtained, i.e.,
\begin{equation}\label{New43}
{\bs s}_{k} = -\frac{{\bs g}_{k}}{\|{\bs g}_{k}\|_2}h_{k},
\end{equation}
which leads to ${\bs a}_{k+1} = {\bs a}_{k} + {\bs s}_{k}$.
If $\|{\bs s}^{C}_{k}\|_2 < h_{k}$ and $\|{\bs s}^{N}_{k}\|_2 > h_{k}$, a dogleg path consisting of two line segments is used to approximate ${\bs s}$ in \eqref{eq2.2.4.4}, i.e.,
\begin{equation}\label{New42}
{\bs s}_{k}(\lambda) = {\bs s}^{C}_{k} + \lambda({\bs s}^{N}_{k} - {\bs s}^{C}_{k}),\quad  0 \leq \lambda \leq 1.
\end{equation}
Obviously, when $\lambda = 0$, ${\bs s}_{k}(\lambda)$ reduces to the steepest descent direction. While $\lambda = 1$, it becomes the Newtonian direction. To exactly obtain $\lambda$ in \eqref{New42}, we will solve the following equation:
\begin{equation*}
\|{\bs s}^{C}_{k} + \lambda({\bs s}^{N}_{k} - {\bs s}^{C}_{k})\|_2 = h_{k}.
\end{equation*}
As a result, we have
\begin{equation*}
{\bs a}_{k+1} = {\bs a}_{k} + {\bs s}_{k}(\lambda) = {\bs a}_{k} + {\bs s}^{C}_{k} + \lambda({\bs s}^{N}_{k} - {\bs s}^{C}_{k}),
\end{equation*}
Otherwise, we choose
\begin{equation}
{\bs s}_{k} = {\bs s}^{N}_{k} = -{\bs G}^{-1}_{k}{\bs g}_{k}. \label{New41}
\end{equation}
In summary, with \eqref{New43}, \eqref{New42} and \eqref{New41}, the solution ${\bs s}_{k}$ in \eqref{eq2.2.4.4} becomes
\begin{equation}
{\bs s}_{k}=\begin{cases}\vspace{0.2cm}
  -\frac{{\bs g}_{k}}{\|{\bs g}_{k}\|_2} h_{k}, &  \textrm{if}\;\; \|{\bs s}^{C}_{k}\|_{2} \geq h_{k},  \quad\quad\quad\quad\quad\quad\quad\quad\quad\quad\quad  (\textrm{I}) \\\vspace{0.2cm}
{\bs s}^{C}_{k} + \lambda({\bs s}^{N}_{k} - {\bs s}^{C}_{k}), & \textrm{if}\;\; \|{\bs s}^{C}_{k} \|_{2} < h_{k}\;\;
 \textrm{and}\;\;  \|{\bs s}^{N}_{k}\|_{2} > h_{k},  \quad\quad\quad (\textrm{II}) \\
 -{\bs G}^{-1}_{k}{\bs g}_{k}, &  \textrm{if}\;\; \|{\bs s}^{C}_{k}\|_{2} < h_{k}\;\;  \textrm{and}\;\;  \|{\bs s}^{N}_{k}\|_{2} \leq h_{k}.  \quad\quad\quad (\textrm{III})
\end{cases}
\end{equation}
\begin{figure}[!ht]
\centering
\subfigure[I]{\includegraphics[width=3.5cm]{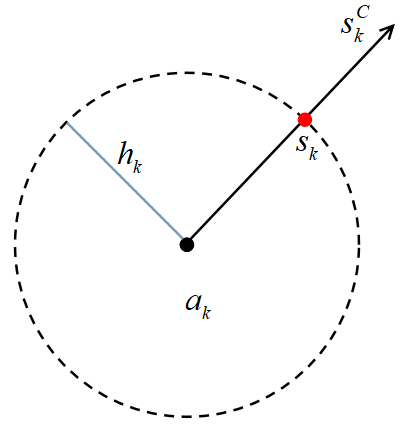}}\;\quad\quad\quad
\subfigure[II]{\includegraphics[width=5cm]{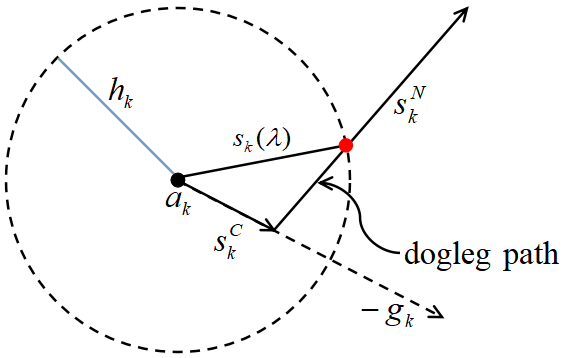}}\;\quad\quad
\subfigure[III]{\includegraphics[width=3.5cm]{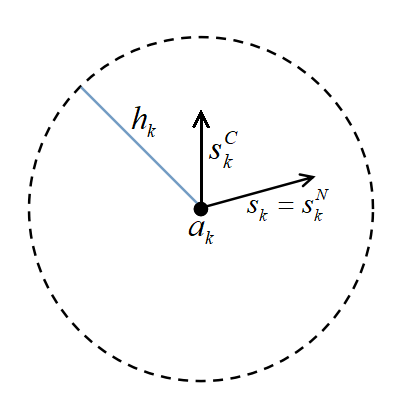}}\;
\caption{Exact trajectory and dogleg approximation.}\label{Example2.1}
\end{figure}
Next, we remark the trust region method for solving nonlinear algebraic system \eqref{eq2.2.4.1}. As mentioned in \cite{sun2006optimization}, the trust region enjoys the desirable global convergence with a local superlinear rate of convergence as follows.
\begin{thm}\label{thm26}  Assume that
\begin{itemize}
 \item[(i)] the function $Q(\textbf{\textit{a}})$ is bounded below on the level set
\begin{equation}
 H := \{\textbf{\textit{a}}\in R^{n}\; :\; Q(\textbf{\textit{a}}) \leq Q(\textbf{\textit{a}}_{0})\}, \quad \forall\, \textbf{\textit{a}}_{0}\in \mathbb R^n,
\end{equation}
and is Lipschitz continuously differentiable in $H;$
\item[(ii)] the Hessian matrixes $G(\textbf{\textit{x}}^{(k)})$ are uniformly bounded in 2-norm, i.e., $\|G(\textbf{\textit{a}}_{k})\| \leq \beta$ for any $ k$ and some  $\beta>0$.
\end{itemize}
If ${\bs g}(\textbf{\textit{a}}_{k}) \neq {\bs 0}$, then
\begin{equation}
\lim_{k \to \infty}\inf\|{\bs g}(\textbf{\textit{a}}_{k})\| = 0.
\end{equation}
Moreover,  if ${\bs g}(\textbf{\textit{a}}^{*}) = {\bs 0}$, and ${\bs G}(\textbf{\textit{a}}^{*})$ is positive definite,  then  the convergence rate of the trust region method is quadratic.
\end{thm}
\begin{rem}\label{remA} \emph{When $k$ is large enough, the trust region method becomes the Newtonian iteration. As a result, it has the same convergence rate as the Newtonian method.} \qed
\end{rem}
\begin{rem} \emph{In practice,  the gradient and Hessian matrices might be appropriately approximated by some numerical means. We refer to Zhang et al. \cite{Hongchaozhang} for such derivative-free methods for \eqref{eq2.2.4.2} with $\bs f$ being twice continuously differentiable, but none of their first-order or second-order derivatives being explicitly
available.} \qed
\end{rem}

\bibliographystyle{plain}
\bibliography{ref}

\end{document}